\documentstyle[12pt,leqno,amstex,amssymb,amsfonts,eucal]{article}
\begin{document}
\oddsidemargin= 20mm
\hoffset=-4pc
\topmargin= 35mm
\textwidth=175mm
\textheight=270mm
\pagestyle{plain}
\newcounter{r}
\newcommand{\Ker}{Ker}
\newcommand{\Gal}{Gal}
\newcommand{\Mat}{Mat}
\newcommand{\sign}{sign}
\newcommand {\Log}{Log}
\newcommand{\Res}{Res}
\newcommand{\Tr}{Tr}
\newcommand{\Nm}{Nm}
\begin{center}
{\large\bf  On the number $\zeta(3)$ }
\end{center}
\begin{center}
{\large\bf L.A.Gutnik}
\end{center}

\vskip.10pt
 {\hfill\sl To the thirtieth anniversary }

\hspace*{8cm}{\hfill\sl of Apery's discovery.}
\vskip.15pt

\begin{center}{\large\bf Table of contents}\end{center}
\vskip.5pt
\S 0. Foreword.
\vskip.5pt\noindent
\S 1. Short proof of Yu.V. Nesterenko expansion.
\vskip.5pt\noindent
\S 2. Introduction. Begin of the proof of Theorem A. 
\vskip.5pt\noindent
\S 3. Transformation the system considered in the Itroduction
in the case $\alpha=1.$
\vskip.5pt\noindent
\S 4. Calculation of the matrix $A^{\ast\ast}_{1,0}(z,\nu).$
\vskip.5pt\noindent
\S 5.  Further properties
 of the matrix $A^{\ast0\ast}_{1,0}(z,\nu).$
\vskip.5pt\noindent
\S 6. Further properties of the functions considered in the Introduction.
\vskip.5pt\noindent
\S 7. End of the proof of theorem A.
{\begin{center}\large\bf\S 0. Foreword . \end{center}}
Let  is given a difference equation
\begin{equation}\label{eq:0c00}
x_{\nu+1}-b_{\nu+1}x_\nu-a_{\nu+1}x_{\nu-1}=0,
\end{equation}
with $\nu\in{\mathbb N}_0.$ We denote by 
$$\{P_\nu(b_0,\,a_1,\,b_1,\,...,\,a_\nu,\,b_\nu)\}_{\nu=-1}^{+\infty}$$
 and
$$\{Q_\nu(b_0,\,a_1,\,b_1,\,...,\,a_\nu,\,b_\nu)\}_{\nu=-1}^{+\infty}$$
the solutions of this equation with initial values
\begin{equation}\label{eq:0d00}
P_{-1}=1,\,Q_{-1}=0,\,P_0(b_0)=b_0,\,Q_{0}(b_0)=1.
\end{equation}
Then
$$ \Bigg\{\frac {P_\nu(b_0,\,a_1,\,b_1,\,...,\,a_\nu,\,b_\nu)}
{Q_\nu(b_0,\,a_1,\,b_1,\,...,\,a_\nu,\,b_\nu)}\Bigg\}_{\nu=0}^{+\infty}$$
\newpage
\pagestyle{headings}
\topmargin= -0.7mm
\textheight=236mm
\vskip.6mm
\markright{\footnotesize\bf L.A.Gutnik, On the number $\zeta(3).$ }
is sequence of convergents of continuous fraction
$$b_0+\frac{a_1\vert}{\vert b_1}+...+\frac{a_\nu\vert}{\vert b_\nu}+...\,\,.$$
Accoding to the famous result of R. Ap\'ery [\ref{r:cd}],  
\begin{equation}\label{eq:0a!}
\zeta(3)=\lim\limits_{\nu\to\infty}\frac{v_\nu}{u_\nu},
\end{equation}
where $\{u_\nu\}_{\nu=1}^{+\infty}$ and $\{v_\nu\}_{\nu=1}^{+\infty}$
 are solutions of difference equation
\begin{equation}\label{eq:0a}
(\nu+1)^3x_{\nu+1}-(34\nu^3+51\nu^2+27\nu+5)x_\nu+\nu^3x_{\nu-1}=0,
\end{equation} 
with initial values $u_0=1,\,u_1=5,\,v_1=0,\,v_1=6.$
The equality (\ref{eq:0a!}) is equivalent to the equality
\begin{equation}\label{eq:0a11}
\zeta(3)=
b_0^\vee+\frac{a_1^\vee\vert}{\vert b_1^\vee}+
\frac{a_2^\vee\vert}{\vert b_2^\vee}+...+
\frac{a_\nu^\vee\vert}{\vert b_\nu^\vee}+...
\end{equation},
with
\begin{equation}\label{eq:g2}
b_0^\vee=0,\,b_1^\vee=5,\,a_1^\vee=6,\,b_{\nu+1}^\vee=
\end{equation}

$$34\nu^3+51\nu^2+27\nu+5,\,a_{\nu+1}^\vee=-\nu^6,$$
where $\nu\in{\mathbb N}.$
Yu.V. Nsterenko in [\ref{r:cf}] has offered the following 
expansion the number $2\zeta(3)$ in continuous fraction:
\begin{equation}\label{eq:N1}
2\zeta(3)=2+\frac{1\vert}{\vert2}+\frac{2\vert}{\vert4}+\frac{1\vert}{\vert3}+
\frac{4\vert}{\vert2}..., 
\end{equation}
with
\begin{equation}\label{eq:N2}
b_0=b_1=a_2=2,\,a_1=1,\,b_2=4,
\end{equation}
\begin{equation}\label{eq:N3}
b_{4k+1}=2k+2,\,a_{4k+1}=k(k+1),\,b_{4k+2}=
\end{equation}
$$2k+4,\,a_{4k+2}=(k+1)(k+2)$$
for $k\in{\mathbb N},$
\begin{equation}\label{eq:N4}
b_{4k+3}=2k+3,\,a_{4k+3}=(k+1)^2,\,b_{4k+4}=
\end{equation}
$$2k+2,\,a_{4k+4}=(k+2)^2$$
for $k\in{\mathbb N}_0.$
The halfs of convergents of continuous fraction (\ref{eq:N1}) compose
a sequence containig convergents of continuous fraction (\ref{eq:0a11}).

Making use of the method developed in our papers
 [\ref{r:caj3}] -- [\ref{r:caj13}], we have found the followng
expansions of the Number $\zeta(3)$ in contiuous fractions :

\bf Teorem A. \it  The followng equalities hold
\newpage
\begin{equation}\label{eq:7.21}
2\zeta(3)=b^{(\ast1)}_0+\frac{a^{(\ast1)}_1\vert}{\vert b^{(\ast1)}_1}+...+
\frac{a^{(\ast1)}_\nu\vert}{\vert b^{(\ast1)}_\nu}+...,
\end{equation}
\begin{equation}\label{eq:7.22}
2\zeta(3)=b^{(\ast2)}_0+\frac{a^{(\ast2)}_1\vert}{\vert b^{(\ast2)}_2}+...+
\frac{a^{(\ast2)}_\nu\vert}{\vert b^{(\ast2)}_\nu}+...,
\end{equation}
with
$$b^{(\ast1)}_0=3,\,a^{(\ast1)}_1=-81,$$
$$a^{(\ast1)}_\nu=-(\nu-1)^3\nu^3(4\nu^2-4\nu-3)^3$$
for $\nu\in[2,+\infty)\cap{\mathbb N},$
$$b^{(\ast1)}_\nu=4(68\nu^6-45\nu^4+12\nu^2-1)$$
for $\nu\in{\mathbb N},$
$$b^{(\ast2)}_0=2,\,a^{(\ast2)}_1=42,$$
$$a^{(\ast2)}_\nu=
-(\nu-1)^3\nu^3(4\nu^2-4\nu-3)((\nu+1)^3-\nu^3)((\nu-1)^3-(\nu-2)^3)
$$
for $\nu\in[2,+\infty)\cap{\mathbb N},$
$$b^{(\ast2)}_\nu=2(102\nu^6-68\nu^4+21\nu^2-3),$$
for $\nu\in{\mathbb N}.$
\rm 
I give short proof of Yu.V. Nesterenko expansion in section 1. I 
prove Theorem A in sections 2 -- 7.

{\begin{center}\large\bf\S 1. Short proof of Yu.V. Nesterenko expansion.
 \end{center}}
Instead of expansion (\ref{eq:N1}) with (\ref{eq:N2}),
 it is more convenient for us to prove the equivalent expansion
\begin{equation}\label{eq:N1a}
\zeta(3)=1+\frac{1\vert}{\vert4}+\frac{4\vert}{\vert4}+\frac{1\vert}{\vert3}+
\frac{4\vert}{\vert2}..., 
\end{equation}
with 
\begin{equation}\label{eq:N2a}
b_0=1,\,a_1=1,\,b_1=a_2=b_2=4.
\end{equation}
Furthermore, to avoid of mishmash in notations, we denote below
$
a_\nu, b_\nu$ for the fraction (\ref{eq:N1a}) by $a_\nu^\wedge,\,b_\nu^\wedge.
$ 
Let $P_{-1}^\vee=1,\,Q_{-1}^\vee=0,$
$$P_\nu^\vee=
P_\nu(b_0^\vee,\,a_1^\vee,\,b_1^\vee,\,...,\,a_\nu^\vee,\,b_\nu^\vee),\,
Q_\nu^\vee=
Q_\nu(b_0^\vee,\,a_1^\vee,\,b_1^\vee,\,...,\,a_\nu^\vee,\,b_\nu^\vee),$$
where values $a_\nu^\vee,\,b_\nu^\vee$ are spcified in (\ref{eq:g2}),
and $\nu\in{\mathbb N}_0.$ Then
\begin{equation}\label{eq:g2a}
Q_0^\vee=1,\,P_0^\vee=b_0^\vee=0,\,
Q_1^\vee=b_1^\vee=5,\,P^\vee_1=
\end{equation}
$$a_1^\vee=6,\,b_2^\vee=117,\,a_2^\vee=-1,P_2^\vee=$$
$$b_2^\vee P^\vee_1+a_2^\vee P^\vee_0=702,Q_2^\vee=$$
$$b_2^\vee Q^\vee_1+a_2^\vee Q^\vee_0=584.$$
Let $P_{-1}^\wedge=1,\,Q_{-1}^\wedge=0,$
$$P_\nu^\wedge=
P_\nu(b_0^\wedge,\,a_1^\wedge,\,b_1^\wedge,\,...,\,a_\nu^\wedge,
\,b_\nu^\wedge),\,
Q_\nu^\wedge=
Q_\nu(b_0^\wedge,\,a_1^\wedge,\,b_1^\wedge,\,...,\,a_\nu^\wedge,
\,b_\nu^\wedge),$$
where $\nu\in{\mathbb N}_0,\,a_\nu^\wedge:=a_\nu,\,b_\nu^\wedge:=b_\nu,$
and values $a_\nu,\,b_\nu$
 are spcified in (\ref{eq:N2a}), (\ref{eq:N3}), and (\ref{eq:N4}). Then,
since $P_{-1}^\wedge=1,\,Q_{-1}^\wedge=0,$ it follows from (\ref{eq:N2a}) that
$$P_0^\wedge=b_0=1,\,Q_0^\wedge=1,\,
P_1^\wedge=b_1^\wedge P_0^\wedge+a_1^\wedge P_{-1}^\wedge=5,\,
Q_1^\wedge=b_1^\wedge Q_0^\wedge+a_1^\wedge Q_{-1}^\wedge=4,$$
\begin{equation}\label{eq:N5}
P_2^\wedge=b_2^\wedge P_1^\wedge+a_2^\wedge P_0^\wedge=24=4P_1^\vee,\,
Q_2^\wedge=b_2^\wedge Q_1+
\end{equation}
$$a_2^\wedge Q_0=20=4Q_1^\vee,\,P_3^\wedge=b_3^\wedge P_2^\wedge+a_3^\wedge P_1^\wedge=77,
\,Q_3^\wedge=b_3^\wedge Q_2^\wedge+$$
$$a_3^\wedge Q_1^\wedge=64,\,P_4^\wedge=b_4^\wedge P_3^\wedge+
a_4^\wedge P_2^\wedge=250,
\,Q_4^\wedge=b_4^\wedge Q_3^\wedge+a_4^\wedge Q_2^\wedge=208,$$
$$P_5^\wedge=b_5^\wedge P_4^\wedge+a_5^\wedge P_3^\wedge=1154,
\,Q_5^\wedge=b_5^\wedge Q_4^\wedge+a_5^\wedge Q_3^\wedge=960,$$
\begin{equation}\label{eq:N6}
P_6^\wedge=b_6^\wedge P_5^\wedge+a_6^\wedge P_4^\wedge=12\times702=12P^\vee_2,
\,Q_6^\wedge=b_6^\wedge Q_5^\wedge+
\end{equation}
$$a_6^\wedge Q_4^\wedge=12\times584=12Q^\vee_2.$$
Let
\begin{equation}\label{eq:N6a}
U_\nu^\wedge=\left(\matrix P_{\nu-1}^\wedge&Q_{\nu-1}^\wedge\\
P_\nu^\wedge&Q_\nu^\wedge\endmatrix\right),
A_\nu^\wedge=\left(\matrix0&1\\a_{1+\nu}^\wedge&
b_{1+\nu}^\wedge\endmatrix\right),
\end{equation}
\begin{equation}\label{eq:N6b}
U_\nu^\vee=\left(\matrix P_{\nu-1}^\vee&Q_{\nu-1}^\vee\\
P_\nu^\vee&Q_\nu^\vee\endmatrix\right),
A_\nu^\vee=\left(\matrix0&1\\a_{1+\nu}^\vee&b_{1+\nu}^\vee\endmatrix\right),
\end{equation}
where $\nu\in{\mathbb N}_0.$
Then
\begin{equation}\label{eq:N7}
U_{\nu}^\wedge=A_{\nu-1}^\wedge U_{\nu-1}^\wedge,
\,U_\nu^\vee=A_{\nu-1}^\vee U_{\nu-1}^\vee
\end{equation}
for $\nu\in{\mathbb N},$
\begin{equation}\label{eq:N8}
U_0^\vee=\left(\matrix1&0\\0&1\endmatrix\right),
U_1^\vee=\left(\matrix0&1\\6&5\endmatrix\right),
U_2^\vee=\left(\matrix 6&5\\702&584\endmatrix\right),
\end{equation}
\begin{equation}\label{eq:N8a}
U_0^\wedge=\left(\matrix1&0\\1&1\endmatrix\right),
U_1^\wedge=\left(\matrix1&1\\5&4\endmatrix\right),
U_2^\wedge=\left(\matrix5&4\\24&20\endmatrix\right),
U_3^\wedge=
\end{equation}
$$\left(\matrix24&20\\77&64\endmatrix\right),U_4^\wedge=
\left(\matrix77&64\\250&208\endmatrix\right),U_5^\wedge=$$
$$\left(\matrix250&208\\1154&960\endmatrix\right),
U_6^\wedge=\left(\matrix1154&960\\8424&7008\endmatrix\right),$$
\begin{equation}\label{eq:N9}
(U_1^\vee)(U_2^\wedge)^{-1}=\frac14\left(\matrix-24&5\\0&1\endmatrix\right),
\end{equation}
\begin{equation}\label{eq:N10}
(U_2^\vee)(U_6^\wedge)^{-1}=\frac1{96}
\left(\matrix-36&5\\0&8\endmatrix\right).
\end{equation}
Let 
\begin{equation}\label{eq:N11}
H_1=\frac14\left(\matrix-24&5\\0&1\endmatrix\right),
\end{equation}
\begin{equation}\label{eq:N12}
H_k=\left(\matrix 12(k+2)(k+1)c(k+1)&-5(k+2)c(k+1)\\
0&-(k-1)^3c(k)\endmatrix\right),
\end{equation}
where  $k\in[2,+\infty)\cap{\mathbb Z}$ and $c(k)=(-2(k-1)^3(k+1)!)^{-1}.$ 
Let $k\in{\mathbb N},\,k\ge2.$
Then, in view of (\ref{eq:N6a}), 
$$A_{4k-6}^\wedge=
\left(\matrix0&1\\a_{4(k-2)+3}^\wedge&b_{4(k-2)+3}^\wedge\endmatrix\right)=
\left(\matrix0&1\\(k-1)^2& 2k-1\endmatrix\right),$$
$$A_{4k-5}^\wedge=\left(\matrix0&1\\a_{4(k-2)+4}^\wedge&b_{4(k-2)+4}^\wedge
\endmatrix\right)=
\left(\matrix0&1\\k^2&2k-2\endmatrix\right),$$
$$A_{4k-4}^\wedge=\left(\matrix0&1\\a_{4(k-1)+1}^\wedge&b_{4(k-1)+1}^\wedge
\endmatrix\right)=
\left(\matrix0&1\\k^2-k&2k\endmatrix\right),$$
$$A_{4k-3}^\wedge=\left(\matrix0&1\\a_{4(k-1)+2}^\wedge&b_{4(k-1)+2}^\wedge
\endmatrix\right)=
\left(\matrix0&1\\k^2+k&2k+2\endmatrix\right),$$
$$A_{4k-5}^\wedge A_{4k-6}^\wedge=\left(\matrix (k-1)^2& 2k-1\\
2(k-1)^3&5k^2-6k+2\endmatrix\right),$$
$$A_{4k-4}^\wedge A_{4k-5}^\wedge A_{4k-6}^\wedge
=\left(\matrix2(k-1)^3&5k^2-6k+2\\
5k(k-1)^3&k(12k^2-15k+5)\endmatrix\right).$$
Let
\begin{equation}\label{eq:N12a}
B_k^\wedge=A_{4k-3}^\wedge A_{4k-4}^\wedge A_{4k-5}^\wedge A_{4k-6}=
\end{equation}
$$\left(\matrix 5k(k-1)^3&k(12k^2-15k+5)\\
12k(k+1)(k-1)^3&k(k+1)(29k^2-36k+12)\endmatrix\right),$$
Then, in view of (\ref{eq:N7}),
\begin{equation}\label{eq:N13}
U_{4k-2}^\wedge=B_k^\wedge U_{4k-6}^\wedge,
\end{equation}
where, as before, $k\in[2,+\infty)\cap{\mathbb Z}.$
Let now $k\in[3,+\infty)\cap{\mathbb Z}.$ Then,
in view of (\ref{eq:N12}),
\begin{equation}\label{eq:N14}
H_{k-1}=\left(\matrix 12(k+1)kc(k)&-5(k+1)c(k)\\
0&-(k-2)^3c(k-1)\endmatrix\right).
\end{equation}

In view of (\ref{eq:g2})
\begin{equation}\label{eq:N15}
\,b_k^\vee=34(k-1)^3+51(k-1)^2+27(k-1)+5=
\end{equation}
$$34k^3-51k^2+27k-5,\,a_k^\vee=-(k-1)^6,$$
where $k\in[3,+\infty)\cap{\mathbb Z}.$ Hence, in view of (\ref{eq:N6b}),
\begin{equation}\label{eq:N16}
A_{k-1}^\vee=\left(\matrix0&1\\
-(k-1)^6&34k^3-51k^2+27k-5\endmatrix\right),
\end{equation}
In view (\ref{eq:N14}) -- (\ref{eq:N16}),
\begin{equation}\label{eq:N17}
A_{k-1}^\vee H_{k-1}=
\left(\matrix0&1\\
-(k-1)^6&34k^3-51k^2+27k-5\endmatrix\right)\times
\end{equation}
$$\left(\matrix 12(k+1)kc(k)&-5(k+1)c(k)\\
0&-(k-2)^3c(k-1)\endmatrix\right)=$$
$$\left(\matrix0&-(k-2)^3c(k-1)\\
-(k-1)^612(k+1)kc(k)&(k-1)^65(k+1)c(k)-b_k^\vee(k-2)^6c(k-1).
\endmatrix\right).$$
In view (\ref{eq:N12a}) and (\ref{eq:N12}),
\begin{equation}\label{eq:N18}
 H_kB_k^\wedge=\left(\matrix 12(k+2)(k+1)c(k+1)&-5(k+2)c(k+1)\\
0&-(k-1)^3c(k)\endmatrix\right)\times
\end{equation}
$$\left(\matrix 5k(k-1)^3&k(12k^2-15k+5)\\
12k(k+1)(k-1)^3&k(k+1)(29k^2-36k+12)\endmatrix\right)=$$
$$\left(\matrix0&(k+2)c(k+1)k(k+1)(-k^2)\\
-c(k)12k(k+1)(k-1)^6&-(k-1)^3c(k)k(k+1)(29k^2-36k+12)
\endmatrix\right).$$
Since
$$-(k+2)(k+1)c(k+1)k^3=-c(k-1)(k-2)^3,$$
$$-(k-1)^3c(k)k(k+1)(29k^2-36k+12)-((k-1)^65(k+1)c(k)=$$
$$-(34k^3-51k^2+27k-5)(k-1)^3(k+1)c(k))=$$
$$-(34k^3-51k^2+27k-5)(k-2)^3)c(k-1)),$$
it follows from (\ref{eq:N15}), (\ref{eq:N17}) and (\ref{eq:N18}) that
\begin{equation}\label{eq:N19}
A_{k-1}^\vee H_{k-1}= H_kB_k^\wedge
\end{equation}
for $k\in[3,+\infty)\cap{\mathbb Z}.$
We prove by induction now the followiong equality
\begin{equation}\label{eq:N20}
U_{k}^\vee=H_kU_{4k-2}^\wedge.
\end{equation}
for any $k\in{\mathbb N}.$
In view of (\ref{eq:N9}) and (\ref{eq:N11}),
the equality (\ref{eq:N20}) holds for $k=1.$
In view of (\ref{eq:N10}) and (\ref{eq:N12}),
the equality (\ref{eq:N20}) holds for $k=2.$
Let $k\in[3,+\infty)\cap{\mathbb Z}$ and (\ref{eq:N20})
holds for $k-1.$ Then, in view of(\ref{eq:N13}), (\ref{eq:N19}) 
 and (\ref{eq:N7}),
$$H_kU_{4k-2}^\wedge=H_kB_kU_{4k-6}^\wedge=
A_{k-1}^\vee H_{k-1}U_{4k-6}^\wedge=A_{k-1}^\vee U_{k-1}^\vee=U_k^\vee.$$
So, the equality (\ref{eq:N20}) holds for any $k\in{\mathbb N}.$
In view of (\ref{eq:N20}), 
\begin{equation}\label{eq:N21}
P_k^\vee=(2(k+1)!)^{-1}P_{4k-2}^\wedge,\,Q_k^\vee=(2(k+1)!)^{-1}Q_{4k-2}^\wedge
\end{equation}
for $k\in[2,+\infty)\cap{\mathbb Z}.$
Since
\begin{equation}\label{eq:N22}
P_\nu^\vee=(\nu!)^3{v_\nu},\,Q_\nu^\vee=(\nu!)^3{u_\nu}
\end{equation}
for $v_\nu$ and $u_\nu$ in (\ref {eq:0a!}) and $\nu\in{\mathbb N}_0,$
it follows from (\ref {eq:N21}) and (\ref {eq:N22}), that
\begin{equation}\label{eq:N23}
P_{4k-2}^\wedge=2(k+1)(k!)^4v_k,\,
Q_{4k-2}^\wedge=2(k+1)(k!)^4u_k
\end{equation}
Let $k\in[2,+\infty)\cap{\mathbb Z}.$ In view of (\ref{eq:N2a}) and
 (\ref{eq:N3}) - (\ref{eq:N4}),
\begin{equation}\label{eq:N24}
Q_{4k-2}^\wedge<Q_{4k-1}^\wedge<Q_{4k-2}^\wedge(k^2+2k+1),
\end{equation}
\begin{equation}\label{eq:N25}
Q_{4k-2}^\wedge<Q_{4k-1}^\wedge<Q_{4k}^\wedge<
\end{equation}
$$Q_{4k-1}^\wedge(k^2+4k+1)<Q_{4k-2}^\wedge(k^2+2k+1)(k^2+4k+1),$$
\begin{equation}\label{eq:N26}
Q_{4k-2}^\wedge<Q_{4k}^\wedge<Q_{4k+1}^\wedge<Q_{4k}^\wedge(k^2+3k+2)<
\end{equation}
$$Q_{4k-2}^\wedge(k^2+2k+1)(k^2+4k+1)(k^2+3k+2),$$
\begin{equation}\label{eq:N27}
\prod\limits_{\kappa=1}^{4k-2} a_\kappa=
4\prod\limits_{\kappa=2}^{k}
a_{4\kappa-5}a_{4\kappa-4}a_{4\kappa-3}a_{4\kappa-2}=
\end{equation}
$$4\prod\limits_{\kappa=2}^{k}(\kappa-1)^2\kappa^2(\kappa-1)\kappa
\kappa(\kappa+1)=2(k!)^8(k+1)/k^3,$$
\begin{equation}\label{eq:N28}
\prod\limits_{\kappa=1}^{4k-1} a_\kappa=
k^2\prod\limits_{\kappa=1}^{4k-2} a_\kappa=2(k!)^8(k+1)/k,
\end{equation}
\begin{equation}\label{eq:N29}
\prod\limits_{\kappa=1}^{4k} a_\kappa=
(k+1)^2\prod\limits_{\kappa=1}^{4k-1} a_\kappa=2(k!)^8(k+1)^3/k,
\end{equation}
\begin{equation}\label{eq:N29a}
\prod\limits_{\kappa=1}^{4k+1} a_\kappa=
k(k+1)\prod\limits_{\kappa=1}^{4k} a_\kappa=2(k!)^8(k+1)^4.
\end{equation}
As it is well known, for any $\varepsilon>0$ there exist
$C_1(\varepsilon)>0$ and $C_2(\varepsilon)>0$ such that
\begin{equation}\label{eq:N30}
C_1(\varepsilon)(1+\sqrt{2})^{4k(1-\varepsilon)}<\vert u_k\vert<
C_2(\varepsilon)(1+\sqrt{2})^{4k(1+\varepsilon)},
\end{equation}     
\begin{equation}\label{eq:N31}
C_1(\varepsilon)(1+\sqrt{2})^{4k(1-\varepsilon)}<\vert v_k\vert<
C_2(\varepsilon)(1+\sqrt{2})^{4k(1+\varepsilon)},
\end{equation}
\begin{equation}\label{eq:N32}
\frac{C_1(\varepsilon)}{(1+\sqrt{2})^{8k(1+\varepsilon)}}<
\bigg\vert\zeta(3)-\frac{v_k}{u_k}\bigg\vert<
\frac{C_2(\varepsilon)}{(1+\sqrt{2})^{8k(1-\varepsilon)}}.
\end{equation}
Terefore, according to (\ref{eq:N23}) -- (\ref{eq:N32}), and 
 well known expression for the difference of two neighboring convergents
of continuous fraction, for any $\varepsilon>0$ there exist
$C_3(\varepsilon)>0$ and $C_4(\varepsilon)>0$ such that
\begin{equation}\label{eq:N33}
\frac{C_3(\varepsilon)}{(1+\sqrt{2})^{8k(1+\varepsilon)}}<
\bigg\vert\zeta(3)-\frac{P_{4k-2}^\wedge}{Q_{4k-2}^\wedge}\bigg\vert
<\frac{C_4(\varepsilon)}{(1+\sqrt{2})^{8k(1-\varepsilon)}},
\end{equation}
\begin{equation}\label{eq:N34}
\frac{C_3(\varepsilon)}{(1+\sqrt{2})^{8k(1+\varepsilon)}}<
\bigg\vert\frac {P_{4k-1+i}^\wedge}{Q_{4k-1+i}^\wedge}-
\frac{P_{4k-2+i}^\wedge}{Q_{4k-2+i}^\wedge}\bigg\vert<
\frac {C_4(\varepsilon))}{(1+\sqrt{2})^{8k(1-\varepsilon)}}
\end{equation}
for $i=0,\,1,\,2,$ and
\begin{equation}\label{eq:N35}
\frac{C_3(\varepsilon)}{(1+\sqrt{2})^{8k(1-\varepsilon)}}<
\bigg\vert \zeta(3)-\frac{P_{4k-2+i}^\wedge}{Q_{4k-2+i}^\wedge}\bigg\vert
<\frac{C_4(\varepsilon)}{(1+\sqrt{2})^{8k(1+\varepsilon)}},
\end{equation}
for $i=0,\,1,\,2,\,3.\,\,\blacksquare .$
{\begin{center}\large\bf\S 2. Introduction. Begin of the proof of Theorem A.
\end{center}}
Let
\begin{equation}\label{eq:1a}
\vert z\vert>1,-3\pi/2<\arg(z)\le\pi/2,\log(z)=\ln(\vert z\vert)+i\arg(z).
\end{equation}
Then $\log(-z)=\log(z)-i\pi,$ if $\Re(z)>0$
 and $\log(z)=\log(-z)-i\pi,$ if $\Re(z)<0.$
Let $\alpha\in{\mathbb N}_0,$  
\begin{equation}\label{eq:b}
f_{\alpha,0,1}^{\ast\vee}(z, \nu):=f_{\alpha,0,1}^{\ast}(z, \nu):=
\sum\limits_{k=0}^{\nu+\alpha}(z)^k\binom{\nu+\alpha}k^2\binom{\nu+k}\nu^2,
\end{equation}
\begin{equation}\label{eq:d}
R(\alpha,t,\nu)=\frac{\prod\limits_{j=1}^\nu(t-j)}
{\prod\limits_{j=0 }^{\nu+\alpha}(t+j)}, 
\end{equation}
\begin{equation}\label{eq:h}
f_{\alpha,0,2}^{\ast}(z, \nu)=
\sum\limits_{t=1}^{+\infty}z^{-t}
\frac{(\nu+\alpha)!^2}{\nu!^2}(R(\alpha,t,\nu))^2,
\end{equation}
\begin{equation}\label{eq:ab2}
f_{\alpha,0,4}^\ast(z,\nu)=-\sum\limits_{t=1}^{+\infty} z^{-t}
\frac{(\nu+\alpha)!^2}{\nu!^2}\left (\frac\partial{\partial t}
(R^2)\right)(\alpha,t,\nu), 
\end{equation}
\begin{equation}\label{eq:ac}
f_{\alpha,0,3}^{\ast}(z,\nu)=(\log(z))f_{\alpha,0,2}^{\ast}(z,\nu)+
f_{\alpha,0,4}^{\ast}(z,\nu),
\end{equation}
\begin{equation}\label{eq:acx}
f_{\alpha,0,k}=\frac{\nu!^2}{(\nu+\alpha)!^2}(z,\nu)
f_{\alpha,0,k}^{\ast}(z,\nu)
\end{equation}
where $k=1,\,2,\,3,\,4,\,\nu\in{\mathbb N}_0.$
Let
\begin{equation}\label{eq:93bb}
\mu=\mu_\alpha(\nu)=(\nu+\alpha)(\nu+1),\,\tau=\tau_\alpha(\nu)
=\nu+\frac{1+\alpha}2,
\end{equation}
\begin{equation}\label{eq:a11}
a_{\alpha,0,1,1}^\vee(z,\nu)=
\frac12(-1+2\alpha-\alpha^2-5\mu+3\alpha\mu-5\mu^2-\alpha\mu^2)+
\end{equation}
$$\frac z2(-4+12\alpha-13\alpha^2+6\alpha^3-\alpha^4)+$$
$$\frac z2\mu(-32+54\alpha-29\alpha^2+5\alpha^3-56\mu+20\alpha\mu),$$
\begin{equation}\label{eq:a12}
a_{\alpha,0,1,2}^\vee(z;\nu)=
-2+3\alpha-\alpha^2-8\mu+\alpha\mu-\alpha^2\mu-4\mu^2+
\end{equation}
$$
z(2-11\alpha+17\alpha^2-10\alpha^3+2\alpha^4-4\mu-11\alpha\mu+
3\alpha^2\mu-20\mu^2),
$$
\begin{equation}\label{eq:a13}
a_{\alpha,0,1,3}^\vee(z;\nu)=-4+
5\alpha-\frac32\alpha^2-\frac12\alpha^3-12\mu-2\alpha\mu+
\end{equation}
$$
z\left(
10-24\alpha+\frac{37}2\alpha^2-\frac{11}2\alpha^3+24\mu-22\alpha\mu
\right),
$$
\begin{equation}\label{eq:a14}
a_{\alpha,0,1,4}^\vee(z;\nu)=(z-1)(6-7\alpha+3\alpha^2+12\mu),
\end{equation}
\begin{equation}\label{eq:b11}
a_{\alpha,0,1,1}^\wedge(z;\nu)=1-\alpha+3\mu+\mu^2+
\end{equation}
$$z(4-8\alpha+5\alpha^2-\alpha^3+24\mu-22\alpha\mu+5\alpha^2\mu+16\mu^2),$$
\begin{equation}\label{eq:b12}
a_{\alpha,0,1,2}^\wedge(z;\nu)=4-2\alpha+8\mu+2\alpha\mu+
\end{equation}
$$z(-4+18\alpha-16\alpha^2+4\alpha^3+16\mu+10\alpha\mu),$$
\begin{equation}\label{eq:b13}
a_{\alpha,0,1,3}^\wedge(z;\nu)=
\end{equation}
$$8-2\alpha+\alpha^2+8\mu+z(-20+28\alpha-5\alpha^2-8\mu),$$
\begin{equation}\label{eq:b14}
a_{\alpha,0,1,4}^\wedge(z;\nu)=-(z-1)(12-2\alpha),
\end{equation}
\begin{equation}\label{eq:a21}
a_{\alpha,0,2,1}^\vee(z;\nu)=
\frac z2(-4+12\alpha-13\alpha^2+6\alpha^3-\alpha^4)+
\end{equation}
$$\frac z2\mu(-32+54\alpha-29\alpha^2+5\alpha^3)+
\frac z2\mu^2(-68+34\alpha-6\alpha^2-24\mu),$$
\begin{equation}\label{eq:a22}
a_{\alpha,0,2,2}^\vee(z;\nu)=
\frac12(-1+2\alpha-\alpha^2-5\mu+3\alpha\mu-5\mu^2-\alpha\mu^2)+
\end{equation}
$$\frac z2
(-10\alpha+21\alpha^2-14\alpha^3+3\alpha^4)+$$
$$\frac z2\mu(-16-20\alpha+17\alpha^2-7\alpha^3-48\mu-28\alpha\mu),$$
\begin{equation}\label{eq:a23}
a_{\alpha,0,2,3}^\vee(z;\nu)=
-2+3\alpha-\alpha^2-8\mu+\alpha\mu-\alpha^2\mu-4\mu^2+
\end{equation}
$$\frac z2
(12-32\alpha+25\alpha^2-5\alpha^3-2\alpha^4)+
z\mu(20-23\alpha-3\alpha^2+4\mu),$$
\begin{equation}\label{eq:a24}
a_{\alpha,0,2,4}^\vee(z;\nu)=
\frac12(-8+10\alpha-3\alpha^2-\alpha^3-24\mu-4\alpha\mu)+
\end{equation}
$$\frac z2( 8-10\alpha+3\alpha^2+\alpha^3+24\mu+4\alpha\mu),$$
\begin{equation}\label{eq:a31}
a_{\alpha,0,3,1}^\vee(z;\nu)=
\frac z2
(-4+12\alpha-13\alpha^2+6\alpha^3-\alpha^4)+
\end{equation}
$$\frac z2\mu(-32+54\alpha-29\alpha^2+5\alpha^3)+$$
$$\frac z2\mu^2(-76+44\alpha-9\alpha^2-\alpha^3-48\mu-4\alpha\mu),$$
\begin{equation}\label{eq:a32}
a_{\alpha,0,3,2}^\vee(z;\nu)=z
(-2+\alpha+4\alpha^2-4\alpha^3+\alpha^4)+
\end{equation}
$$z\mu(-16-\alpha+7\alpha^2-3\alpha^3-\alpha^4-34\mu-
17\alpha\mu-7\alpha^2\mu-12\mu^2),$$
\begin{equation}\label{eq:a33}
a_{\alpha,0,3,3}^\vee(z;\nu)=
\frac12(-1+2\alpha-\alpha^2-5\mu+3\alpha\mu-5\mu^2-\alpha\mu^2)+
\end{equation}
$$\frac z2
(4-16\alpha+15\alpha^2-4\alpha^3-\alpha^5)+$$
$$\frac z2\mu(16-42\alpha+\alpha^2-9\alpha^3+8\mu-20\alpha\mu),$$
\begin{equation}\label{eq:a34}
a_{\alpha,0,3,4}^\vee(z;\nu)=
-2+3\alpha-\alpha^2-8\mu+\alpha\mu-\alpha^2\mu-4\mu^2+
\end{equation}
$$z(2-3\alpha+\alpha^2+8\mu-\alpha\mu+\alpha^2\mu+4\mu^2),$$
\begin{equation}\label{eq:a41}
a_{\alpha,0,4,1}^\vee(z;\nu)=
\frac z2
(-4+12\alpha-13\alpha^2+6\alpha^3-\alpha^4)+
\end{equation}
$$\frac z2\mu(-32+54\alpha-29\alpha^2+5\alpha^3-80\mu+50\alpha\mu-
11\alpha^2\mu-\alpha^3\mu)+$$
$$z\mu^3(-32-\alpha-\alpha^2-4\mu),$$
\begin{equation}\label{eq:a42}
a_{\alpha,0,4,2}^\vee(z;\nu)=\frac z2
(-8+14\alpha-5\alpha^2-2\alpha^3+\alpha^4)+
\end{equation}
$$\frac z2\mu(-56+32\alpha+\alpha^2-5\alpha^3-2\alpha^4)+$$
$$\frac z2\mu^2(-112-26\alpha-15\alpha^2-5\alpha^3)+
z\mu^3(-28-10\alpha),$$
\begin{equation}\label{eq:a43}
a_{\alpha,0,4,3}^\vee(z;\nu)=\frac z2
(-4+5\alpha^2-2\alpha^3-\alpha^5)+
\end{equation}
$$\frac z2\mu(-24-22\alpha-3\alpha^2-9\alpha^3-4\alpha^4-36\mu-
42\alpha\mu-18\alpha^2\mu-8\mu^2),$$
\begin{equation}\label{eq:a44}
a_{\alpha,0,4,4}^\vee(z;\nu)=
\frac12(-1+2\alpha-\alpha^2-5\mu+3\alpha\mu-5\mu^2-\alpha\mu^2)+
\end{equation}
$$\frac z2(-2\alpha+\alpha^2-\alpha^5)+
\frac z2\mu(-8\alpha-5\alpha^2-5\alpha^3-4\alpha\mu),$$
\begin{equation}\label{eq:b21}
a_{\alpha,0,2,1}^\wedge(z;\nu)=z(4-8\alpha+5\alpha^2-\alpha^3)+
\end{equation}
$$z\mu(24-22\alpha+5\alpha^2+28\mu-2\alpha\mu),$$
\begin{equation}\label{eq:b22}
a_{\alpha,0,2,2}^\wedge(z;\nu)=1-\alpha+3\mu+\mu^2+
\end{equation}
$$z(10\alpha-11\alpha^2+3\alpha^3)+
z\mu(16+16\alpha+\alpha^2)+16z\mu^2,$$
\begin{equation}\label{eq:b23}
a_{\alpha,0,2,3}^\wedge(z;\nu)=4-2\alpha+8\mu+2\alpha\mu+
\end{equation}
$$z(-12+20\alpha-5\alpha^2+2\alpha^3)+z\mu(-16+14\alpha),$$
\begin{equation}\label{eq:b24}
a_{\alpha,0,2,4}^\wedge(z;\nu)=(8-2\alpha+\alpha^2+8\mu)(1-z),
\end{equation}
\begin{equation}\label{eq:b31}
a_{\alpha,0,3,1}^\wedge(z;\nu)=z(4-8\alpha+5\alpha^2-\alpha^3)+
\end{equation}
$$z\mu(24-22\alpha+5\alpha^2+36\mu-4\alpha\mu+\alpha^2\mu+8\mu^2),$$
\begin{equation}\label{eq:b32}
a_{\alpha,0,3,2}^\wedge(z;\nu)=z(4+2\alpha-6\alpha^2+2\alpha^3)+
\end{equation}
$$z\mu(24+14\alpha+2\alpha^3)+z\mu^2(28+14\alpha),$$
\begin{equation}\label{eq:b33}
a_{\alpha,0,3,3}^\wedge(z;\nu)=1-\alpha+3\mu+\mu^2+
\end{equation}
$$z(-4+12\alpha-3\alpha^2+\alpha^3+\alpha^4-8\mu+18\mu\alpha+7\mu\alpha^2),$$
\begin{equation}\label{eq:b34}
a_{\alpha,0,3,4}^\wedge(z;\nu)=(4-2\alpha+8\mu+2\alpha\mu)(1-z),
\end{equation}
\begin{equation}\label{eq:b41}
a_{\alpha,0,4,1}^\wedge(z;\nu)=z(4-8\alpha+5\alpha^2-\alpha^3)+
\end{equation}
$$
z\mu(24-22\alpha+5\alpha^2+40\mu-6\alpha\mu+\alpha^2\mu)+z\mu^3(16+2\alpha),
$$
\begin{equation}\label{eq:b42}
a_{\alpha,0,4,2}^\wedge(z;\nu)=z(8-6\alpha-\alpha^2+\alpha^3)+
\end{equation}
$$
z\mu(40+4\alpha+\alpha^2+2\alpha^3)+z\mu^2(48+22\alpha+5\alpha^2)+8z\mu^3,$$
\begin{equation}\label{eq:b43}
a_{\alpha,0,4,3}^\wedge(z;\nu)=z(4+4\alpha-\alpha^2+\alpha^3+\alpha^4)+
\end{equation}
$$z\mu(16+22\alpha+11\alpha^2+4\alpha^3)+z\mu^2(12+10\alpha)$$
\begin{equation}\label{eq:b44}
a_{\alpha,0,4,4}^\wedge(z;\nu)=1-\alpha+3\mu+\mu^2+
\end{equation}
$$2z\alpha+z\alpha^2+z\alpha^3+z\alpha^4+4z\mu\alpha+3z\mu\alpha^2,$$

\begin{equation}\label{eq:abik}
a_{\alpha,0,i,k}^\ast(z;\nu)=a_{\alpha,0,i,k}^\vee(z;\nu)+
\tau a_{\alpha,0,i,k}^\wedge(z;\nu),
\end{equation}
where $i=1,\,...,\,4,\,k=1,\,...,\,4.$
We denote by
$$A_{\alpha,0}^\ast(z;\nu),\,A_{\alpha,0}^\vee(z;\nu)(z;\nu),
\,A_{\alpha,0}^\wedge(z;\nu)(z;\nu)$$
the $4\times4$-matrix,  such that its
element in $i$-th row and $k$-th column is equal respectively to  
$$a_{\alpha,0,i,k}^\ast(z;\nu),\,a_{\alpha,0,i,k}^\vee(z;\nu),\,
a_{\alpha,0,i,k}^\wedge(z;\nu)$$
where $i=1,\,...,\,4,\,k=1,\,...,\,4.$
Then
\begin{equation}\label{eq:15c}
A_{\alpha,0}^\ast(z;\nu)=
A_{\alpha,0}^\vee(z;\nu)(z;\nu)+\tau A_{\alpha,0}^\wedge(z;\nu)(z;\nu).
\end{equation}
Clearly,
\begin{equation}\label{eq:15d}
A_{\alpha,0}^\ast(z;\nu)=
A_{\alpha,0}^\ast(z;\nu)(1;\nu)+(z-1)V_{\alpha,0}^\ast(\nu),
\end{equation}
where the matrix $V_{\alpha,0}^\ast(\nu)$ does not depend from $z.$
Let
\begin{equation}\label{eq:15e}
X_{\alpha,0,k}(z;\nu)=
\left(\matrix
f_{\alpha,0,k}(z,\nu)\\
\delta f_{\alpha,0,k}(z,\nu)\\
\delta^2 f_{\alpha,0,k}(z,\nu)\\
\delta^3 f_{\alpha,0,k}(z,\nu)
\endmatrix\right)
,\,X_{\alpha,0,k}^\ast(z;\nu)=\frac{((\nu+\alpha)!)^2}{(\nu!)^2}
X_{\alpha,0,k}(z;\nu)
\end{equation}
for
$k=1,\,2,\,3,\,\vert z\vert>1,\,\nu\in{\mathbb N}_0.$
Let further
\begin{equation}\label{eq:93be0}
X_{\alpha,0,k}(z;-\nu-1-\alpha)=X_{\alpha,0,k}(z;\nu),
\end{equation}
where $\nu\in{\mathbb N}_0.$
The following results are obtained in [\ref{r:caj11}] -- [\ref{r:caj13}].

\bf Theorem 1. \it The column $X_{\alpha,0,k}(z;\nu)$ satisfies
to the equation
\begin{equation}\label{eq:93ca}
\nu^5X_{\alpha,0,k}(z;\nu-1)=
A_{\alpha,0}^\ast(z;\nu)X_{\alpha,0,k}(z;\nu),
\end{equation}  
for
$\nu\in M_\alpha^\ast=(-\infty,-1-\alpha]\cup[1,+\infty))\cap{\mathbb Z},
\,k=1,\,2,\,3,\,\vert z\vert>1;$ moreover,
 the matrix $A_{\alpha,0}^\ast(z;\nu)$ has the following property:
\begin{equation}\label{eq:93cc}
-\nu^5(\nu+\alpha)^5E_{4}=
A_{\alpha,0}^\ast(z;-\nu-\alpha)A_{\alpha,0}^\ast(z;\nu),
\end{equation}
where $E_4$ is the $4\times4$ unit matrix,
$z\in{\mathbb C},\,\nu\in{\mathbb C}.$
\rm
Let us consider the row
\begin{equation}\label{eq:113bd1}
R_{\alpha,0}(\nu)=
(r_{\alpha,0,1}(\nu),\,r_{\alpha,0,2}(\nu),\,
r_{\alpha,0,3}(\nu),\,r_{\alpha,0,4}(\nu)),
\end{equation}
where
\begin{equation}\label{eq:93bd}
r_{\alpha,0,1}(\nu)=\mu_\alpha(\nu)^2,\,r_{\alpha,0,2}(\nu)=
-2(1-\alpha)\times
\end{equation}
$$\mu_\alpha(\nu),\,r_{\alpha,0,3}(\nu)=
(1-\alpha)^2-2\mu_\alpha(\nu),\,r_{\alpha,0,4}(0,\nu)=2(1-\alpha).$$
The following Theorem is proved in [\ref{r:caj13}] (Lemma 11.3.1).

\bf Theorem 2. \it The row $R_{\alpha,0}(\nu)$
has the following property:
\begin{equation}\label{eq:15f}
R_{\alpha,0}(\nu-1)A_{\alpha,0}^\ast(1;\nu)=\nu^5R_{\alpha,0}(\nu),
\end{equation}
where $\nu\in{\mathbb C}.$\rm
{\begin{center}\large\bf\S 3. Transformations of the system considered in 
the Introduction
in the case $\alpha=1.$\end{center}}
In view of (\ref{eq:93bb}), (\ref{eq:93bd})
\begin{equation}\label{eq:15.1.a}
\tau=\tau_1(\nu)=\nu+1,\,\mu=\mu_1(\nu)=(\nu+1)^2,\,r_{1,0,1}(\nu)=
\end{equation}
$$\mu_1(\nu)^2=(\nu+1)^4=\tau^4,\,r_{1,0,2}(\nu)=0,\,r_{1,0,3}(\nu)=$$
$$-2\mu_1(\nu)=-2(\nu+1)^2,\,r_{\alpha,0,4}(0,\nu)=0.$$
Let $E_4$ denotes $4\times4$-unit matrix, and let  $C^{(\kappa)}(\nu)$
 with $\kappa=0,\,1$
is result of replacement of $2\kappa+1$-th row of the matrix $E_4$ 
 by the row in (\ref{eq:113bd1})
with $\alpha=1.$ Let further $D^{(\kappa)}(\nu)$ with $\kappa=0,\,1$ denotes
the adjoint matrix to the matrix $C^{(\kappa)}(\nu).$
 Then
\begin{equation}\label{eq:15.1.b}
C^{(0)}(\nu)=
 \end{equation}
$$\left(\matrix
 r_{1,0,1}(\nu)&r_{1,0,2}(\nu)&r_{1,0,3}(\nu)&r_{1,0,4}(\nu)\\
      0        &     1        &     0        &     0        \\
      0        &     0        &     1        &     0        \\
      0        &     0        &     0        &     1
 \endmatrix\right),\,C^{(1)}(\nu)=$$
 $$
 \left(\matrix
      1        &     0        &    0        &      0        \\
      0        &     1        &    0        &      0        \\
 r_{1,0,1}(\nu)&r_{1,0,2}(\nu)&r_{1,0,3}(\nu)&r_{1,0,4}(\nu)\\
      0        &     0        &     0        &     1
 \endmatrix\right),$$
\begin{equation}\label{eq:15.1.c}
D^{(0)}(\nu)=
\end{equation}
$$\left(\matrix
    1&-r_{1,0,2}(\nu)&-r_{1,0,3}(\nu)&-r_{1,0,4}(\nu)\\
    0& r_{1,0,1}(\nu)&      0        &      0        \\
    0&     0         & r_{1,0,1}(\nu)&      0        \\
    0&     0         &      0        & r_{1,0,1}(\nu) 
\endmatrix\right),\,D^{(1)}(\nu)=$$
$$\left(\matrix
 r_{1,0,3}(\nu)&      0        &0&      0        \\
      0        & r_{1,0,3}(\nu)&0&      0        \\
-r_{1,0,1}(\nu)&-r_{1,0,2}(\nu)&1&-r_{1,0,4}(\nu)\\
      0        &      0        &0& r_{1,0,3}(\nu) 
\endmatrix\right).$$
Clearly,
\begin{equation}\label{eq:15.1.d}
C^{(\kappa)}(\nu)D^{(\kappa)}(\nu)=
(-2)^\kappa(\mu_1(\nu))^{2-\kappa}E_4,\,C^{(\kappa)}(-\nu-2)=
\end{equation}
$$C^{(\kappa)}(\nu),\,D^{(\kappa)}(-\nu-2)=D^{(\kappa)}(\nu),$$
where $\kappa=0,\,1.$ Let
\begin{equation}\label{eq:15.1.e}
A^{\ast\kappa\ast}_{1,0}(z,\nu)=
C^{(\kappa)}(\nu-1)A^{\ast}_{1,0}(z,\nu)D^{(\kappa)}(\nu).
\end{equation}
Then
\begin{equation}\label{eq:15.1.f}
A^{\ast\kappa\ast}_{1,0}(z,-\nu-1)=
\end{equation}
$$C^{(\kappa)}(-\nu-2)A^{\ast}_{1,0}(z,-\nu-1)D^{(\kappa)}(-\nu-1)=$$
$$C^{(\kappa)}(\nu)A^{\ast}_{1,0}(z,-\nu-1)D^{(\kappa)}(\nu-1),$$
and, in view of (\ref{eq:15.1.d}), (\ref{eq:93cc}), (\ref{eq:15.1.e}),
\begin{equation}\label{eq:15.1.g}
A^{\ast\kappa\ast}_{1,0}(z,-\nu-1)A^{\ast\kappa\ast}_{1,0}(z,\nu)=
\end{equation}
$$C^{(\kappa)}(\nu)A^{\ast}_{1,0}(z,-\nu-1)D^{(\kappa)}(\nu-1)
C^{(k)}(\nu-1)A^{\ast}_{1,0}(z,\nu)D(\nu)=$$
$$-4^\kappa(\mu_1(\nu)\mu_1(\nu-1))^{2-\kappa}(\nu(\nu+1))^5E_4,$$
where $\kappa=0,\,1.$ 
Let
\begin{equation}\label{eq:15.1.h}
Y^{(\kappa)}_{1,0,k}(z;\nu)=C^{(\kappa)}(\nu)X_{1,0,k}(z;\nu),
\end{equation}
where
$\kappa=0,\,1,\,k=1,\,2,\,3,\,\vert z\vert>1,\,
\nu\in M_1^{\ast\ast\ast}=((-\infty,-2]\cup[0,+\infty))\cap{\mathbb Z}.$
Then, in view of (\ref{eq:93be0}), (\ref{eq:15.1.d}), (\ref{eq:93ca}),
\begin{equation}\label{eq:15.1.i}
Y^{(\kappa)}_{1,0,k}(z;-\nu-2)=Y^{(\kappa)}_{1,0,k}(z;\nu),
\end{equation}
\begin{equation}\label{eq:15.1.aj}
A^{\ast\kappa\ast}_{1,0}(z,\nu)Y^{(\kappa)}_{1,0,k}(z;\nu)=
\end{equation}
$$C^{(\kappa)}(\nu-1)A^{\ast}_{1,0}(z,\nu)D^{(\kappa)}(\nu)
C^{(\kappa)}(\nu)X_{1,0,k}(z;\nu)=$$
$$(-1)^\kappa(\kappa+1)\mu_1(\nu)^{2-\kappa}
C^{(\kappa)}(\nu-1)A^{\ast}_{1,0}(z,\nu)X_{1,0,k}(z;\nu)=$$
$$
(-1)^\kappa(\kappa+1)\mu_1(\nu)^{2-\kappa}\nu^5
C^{(\kappa)}(\nu-1)X_{1,0,k}(z;\nu-1)=
$$
$$
(-1)^\kappa(\kappa+1)\mu_1(\nu)^{2-\kappa}\nu^5\nu^5
Y^{(\kappa)}_{1,0,k}(z;\nu-1),$$
where
$\kappa=0,\,1,\,k=1,\,2,\,3,\,\vert z\vert>1,\,
\nu\in M_1^{\ast}=((-\infty,-2]\cup[1,+\infty))\cap{\mathbb Z}.$
Replacing in the equality (\ref{eq:15.1.aj})
$$\nu\in M_1^{\ast}=((-\infty,-2]\cup[1,+\infty))\cap{\mathbb Z}$$
 by
$$\nu:=-\nu-2\in M_1^{\ast\ast}=((-\infty,-3]\cup[0,+\infty))\cap{\mathbb Z},$$
and taking in account (\ref{eq:15.1.i}) we obtain 
the equality
\begin{equation}\label{eq:15.1.ba}
-A^{\ast\kappa\ast}_{1,0}(z,-\nu-2)Y^{(\kappa)}_{1,0,k}(z;\nu)=
(-1)^\kappa(\kappa+1)\mu_1(\nu)^{2-\kappa}(\nu+2)^5
Y^{(\kappa)}_{1,0,k}(z;\nu+1),
\end{equation}
where $\kappa=0,\,1,\,k=1,\,2,\,3,\,\vert z\vert>1,\,
\nu\in M_1^{\ast\ast}=((-\infty,-3]\cup[0,+\infty))\cap{\mathbb Z}.$
{\begin{center}\large\bf
\S 4. Calculation of the matrix $A^{\ast\kappa\ast}_{1,0}(z,\nu)$
for $\kappa=0,\,1.$
\end{center}}
Let
\begin{equation}\label{eq:15.2.a}
V^{\ast\kappa\ast}_{1,0}(\nu)=C^{(\kappa)}(\nu-1)V^{\ast}_{1,0}(\nu)
D^{(\kappa)}(\nu).
\end{equation}
Then, in view of $(\ref{eq:15d}),$
 \begin{equation}\label{eq:15.2.b} 
A_{\alpha,0}^{\ast\kappa\ast}(z;\nu)=
A_{\alpha,0}^{\ast\kappa\ast}(1;\nu)(1;\nu)+(z-1)
V_{\alpha,0}^{\ast\kappa\ast}(\nu),
\end{equation}
where the matrix $V_{\alpha,0}^{\ast\ast}(\nu)$ does not depend from $z.$ 
We note that the $2\kappa+1$-th row of the matrix 
$$C^{(\kappa)}(\nu-1)A^{\ast}_{1,0}(1,\nu)$$
coincides with the row $$R_{1,0}^\ast(\nu-1)A^{\ast\ast}_{1,0}(z,\nu)$$
and, according to the Theorem 2 coincides
 with the row $\nu^5R_{1,0}^\ast(\nu),$ i.e. with
the $2\kappa+1$-th row of the matrix $\nu^5C^{(\kappa)}(\nu).$
 Therefore, in view 
of (\ref{eq:15.1.d}), the $(2\kappa+1)$-th row of the matrix
 $A^{\ast\kappa\ast}_{1,0}(1,\nu)$
 is equal to
 \begin{equation}\label{eq:v1}
 (-2)^\kappa(\mu_1(\nu)^{2-\kappa})\nu^5\bar e_{4,2\kappa+1},
 \end{equation}
   where $\bar e_{4,l}$
denotes the $l$-th row of the matrix $E_4$
for $l=1,\,2,\,3,\,4.$   We note further that the second,
 $3-2\kappa$-th and fourth row of the
  matrix $C^{(\kappa)}(\nu-1)A^{\ast}_{1,0}(1,\nu)$
  coincides with respectively  the second, $3-2\kappa$-th and fourth row
 of the matrix $A^{\ast}_{1,0}(1,\nu).$
In view of the equalities (148) -- (151) in [\ref{r:caj13}], since $\alpha=1$
and $\mu=\mu_1(\nu)=(\nu+1)^2,\,\tau=\tau_1(\nu)=\nu+1$ now,
it follows that
\begin{equation}\label{eq:15.2.c}
\nu^5=p_5(\mu)+\tau q_5(\mu),
\end{equation}
where
\begin{equation}\label{eq:15.2.d}
p_5(\mu)=-5\mu^2-10\mu-1,\,q_5(\mu)=\mu^2+10\mu+5.
\end{equation}
Let $A_{1,0}^{\vee\kappa\vee}(1;\nu)$  denotes the $4\times4$-matrix
with second, $(3-2\kappa)$-th and fourth
row equal respectively to the second, $(3-2\kappa)$-th and fourth
row of the matrix $C(\nu-1)A^{\vee}_{1,0}(1,\nu)D(\nu)$ and
 with $(2\kappa+1)$-th row equal to 
$$(-2)^\kappa(\mu_1(\nu)^{2-\kappa})p_5(\mu)\bar e_{4,2\kappa+1}.$$
Let further $A_{1,0}^{\wedge\kappa\wedge}(1;\nu)$ denotes the $4\times4$-matrix
with second, $(3-2\kappa)$-th and fourth
row equal respectively to the second, $(3-2\kappa)$-th and fourth
row of the matrix $C(\nu-1)A^{\wedge}_{1,0}(1,\nu)D(\nu)$ and
 with $(2\kappa+1)$-th row equal to 
$$(-2)^\kappa(\mu_1(\nu)^{2-\kappa})q_5(\mu)\bar e_{4,2\kappa+1}.$$
Then, in view of (\ref{eq:15c}), (\ref{eq:v1}) and (\ref{eq:15.2.c}),
\begin{equation}\label{eq:15.2.e}
A_{1,0}^{\ast\kappa\ast}(1;\nu)=
A_{1,0}^{\vee\kappa\vee}(1;\nu)+\tau A_{1,0}^{\wedge\kappa\wedge}(1;\nu),
\end{equation}
where $\kappa=0,\,1.$
We denote by
$a_{1,0,i,j}^{\vee\kappa\vee}(1;\nu),\,
a_{1,0,i,j}^{\wedge\kappa\wedge}(1;\nu)$ and
$a_{1,0,i,j}^{\ast\kappa\ast}(1;\nu),$
where $\kappa=0,1;\,i,j=1,2,3,4,$
the expressions, which stand on intersection of $i$-th row and $j$-th column
in the matrices respectively
$A_{1,0}^{\ast\kappa\ast}(1;\nu),\,
A_{1,0}^{\vee\kappa\vee}(1;\nu)$ and $A_{1,0}^{\wedge\kappa\wedge}(1;\nu).$
Then, in view of (\ref{eq:v1}),
\begin{equation}\label{eq:15.2.a11}
a_{1,0,2\kappa+1,2\kappa+1}^{\ast\kappa\ast}(1;\nu)=
\end{equation}
$$(-2)^\kappa\tau^{4-2\kappa}(\tau-1)^5,$$
\begin{equation}\label{eq:15.2.vw1,234}
a_{1,0,2\kappa+1,k}^{\ast\kappa\ast}(1;\nu)=0,
\end{equation}
where $\kappa=0,\,1,\,k=3-2\kappa,\,2,\,4.$
In view of (\ref{eq:a11}), 
\begin{equation}\label{eq:av11}
a_{1,0,1,1}^\vee(1;\nu)=-2\mu-21\mu^2.
\end{equation}
In view of (\ref{eq:b11}), 
\begin{equation}\label{eq:bw11}
a_{1,0,1,1}^\wedge(1;\nu)=10\mu+17\mu^2.
\end{equation}
In view of (\ref{eq:av21}), (\ref{eq:bw21}), (\ref{eq:15c})
 and (\ref{eq:15.2.e}),
\begin{equation}\label{eq:as11}
a_{1,0,1,1}^\ast(1;\nu)=
\end{equation}
$$17\tau^5-21\tau^4+10\tau^3-2\tau^2)=\tau^2(17\tau^3-21\tau^2+10\tau-2).$$
In view of (\ref{eq:a21}), (\ref{eq:15.1.b}), (\ref{eq:15.1.c}),
\begin{equation}\label{eq:av21}
a_{1,0,2,1}^{\vee0\vee}(1;\nu)=a_{1,0,2,1}^\vee(1;\nu)=-\mu-20\mu^2-12\mu^3.
\end{equation}
In view of (\ref{eq:b21}), (\ref{eq:15.1.b}), (\ref{eq:15.1.c}),
\begin{equation}\label{eq:bw21}
a_{1,0,2,1}^{\wedge0\wedge}(1;\nu)=a_{1,0,2,1}^\wedge(1;\nu)=7\mu+26\mu^2.
\end{equation}
In view of (\ref{eq:av21}), (\ref{eq:bw21}), (\ref{eq:15c})
 and (\ref{eq:15.2.e}),
\begin{equation}\label{eq:as21}
a_{1,0,2,1}^{\ast0\ast}(1;\nu)=a_{1,0,2,1}^\ast(1;\nu)=
\end{equation}
$$-(12\tau^6-26\tau^5+20\tau^4-7\tau^3+\tau^2)=$$
$$-\tau^2(\tau-1)(12\tau^3-14\tau^2+6\tau-1)=$$
$$-\tau^2(\tau-1)(2\tau-1)(6\tau^2-4\tau+1).$$
In view of (\ref{eq:a31}), (\ref{eq:15.1.b}), (\ref{eq:15.1.c}),
\begin{equation}\label{eq:av31}
a_{1,0,3,1}^{\vee0\vee}(1;\nu)=a_{1,0,3,1}^\vee(1;\nu)=-\mu-21\mu^2-26\mu^3.
\end{equation}
In view of (\ref{eq:b31}), (\ref{eq:15.1.b}), (\ref{eq:15.1.c}),
\begin{equation}\label{eq:bw31}
a_{1,0,3,1}^{\wedge0\wedge}(1;\nu)=a_{\alpha,0,3,1}^\wedge(1;\nu)=
\end{equation}
$$7\mu+33\mu^2+8\mu^3.$$
In view of (\ref{eq:av31}), (\ref{eq:bw31}),  (\ref{eq:15c})
 and (\ref{eq:15.2.e}),
\begin{equation}\label{eq:as31}
a_{1,0,3,1}^{\ast0\ast}(1;\nu)=a_{1,0,3,1}^\ast(1;\nu)=
\end{equation}
$$8\tau^7-26\tau^6+33\tau^5-21\tau^4+$$
$$\tau^3-\tau^2=\tau^2(\tau-1)(8\tau^4-18\tau^3+15\tau^2-6\tau+1)=$$
$$\tau^2(\tau-1)^2(8\tau^3-10\tau^2+5\tau-1).$$
$$\tau^2(\tau-1)^2(2\tau-1)(4\tau^2-3\tau+1).$$
In view of (\ref{eq:a41}), (\ref{eq:15.1.b}), (\ref{eq:15.1.c}),
\begin{equation}\label{eq:av41}
a_{1,0,4,1}^{\vee0\vee}(1;\nu)=a_{1,0,4,1}^\vee(1;\nu)=
-\mu-21\mu^2-34\mu^3-4\mu^4.
\end{equation}
In view of (\ref{eq:b41}), (\ref{eq:15.1.b}), (\ref{eq:15.1.c}),
\begin{equation}\label{eq:bw41}
a_{\alpha,0,4,1}^{\wedge0\wedge}(1;\nu)=a_{\alpha,0,4,1}^{\wedge}(1;\nu)=
7\mu+35\mu^2+18\mu^3.
\end{equation}
In view of (\ref{eq:av41}), (\ref{eq:bw41}), 
 (\ref{eq:15c}) and (\ref{eq:15.2.e}),
\begin{equation}\label{eq:as41}
a_{1,0,4,1}^{\ast0\ast}(1;\nu)=a_{1,0,4,1}^\ast(1;\nu)=
\end{equation}
$$-4\tau^8+18\tau^7-34\tau^6+35\tau^5-21\tau^4+$$
$$
7\tau^3-\tau^2=-\tau^2(\tau-1)(4\tau^5-14\tau^4+20\tau^3-15\tau^2+6\tau-1)=
$$
$$-\tau^2(\tau-1)^2(4\tau^4-10\tau^3+10\tau^2-5\tau+1)=$$
$$-\tau^2(\tau-1)^3(4\tau^3-6\tau^2+4\tau-1)=$$
$$-\tau^2(\tau-1)^3(2\tau-1)(2\tau^2-2\tau+1).$$
In view of (\ref{eq:a13}), (\ref{eq:15.1.b}), (\ref{eq:15.1.c}),
\begin{equation}\label{eq:a1v13}
a_{1,0,1,3}^{\vee1\vee}(1;\nu)=a_{1,0,1,3}^\vee(1;\nu)=-2-12\mu.
\end{equation}
In view of (\ref{eq:b13}), (\ref{eq:15.1.b}), (\ref{eq:15.1.c}),
\begin{equation}\label{eq:b1w13}
a_{1,0,1,3}^{\wedge1\wedge}(1;\nu)=a_{1,0,1,3}^\wedge(1;\nu)=10,
\end{equation}
In view of (\ref{eq:a1v13}), (\ref{eq:b1w13}), (\ref{eq:15c})
 and (\ref{eq:15.2.e}),
\begin{equation}\label{eq:a1s13}
a_{1,0,1,3}^{\ast1\ast}(1;\nu)=a_{1,0,1,3}^\ast(1;\nu)=
\end{equation}
$$-12\tau^2+10\tau-2=-2(2\tau-1)(3\tau-1).$$
In view of (\ref{eq:a23}), (\ref{eq:15.1.b}), (\ref{eq:15.1.c}),
\begin{equation}\label{eq:a1v23}
a_{1,0,2,3}^{\vee1\vee}(1;\nu)=a_{1,0,2,3}^\vee(1;\nu)=-1-14\mu.
\end{equation}
In view of (\ref{eq:b23}), (\ref{eq:15.1.b}), (\ref{eq:15.1.c}),
\begin{equation}\label{eq:b1w23}
a_{1,0,2,3}^{\wedge1\wedge}(1;\nu)=a_{\alpha,0,2,3}^\wedge(1;\nu)=
7+8\mu,
\end{equation}
In view of (\ref{eq:a1v23}), (\ref{eq:b1w23}), (\ref{eq:15c})
 and (\ref{eq:15.2.e}),
\begin{equation}\label{eq:a1s23}
a_{1,0,2,3}^{\ast1\ast}(1;\nu)=a_{1,0,2,3}^\ast(1;\nu)=
\end{equation}
$$8\tau^3-14\tau^2+7\tau-1=$$
$$(\tau-1)(8\tau^2-6\tau+1)=(\tau-1)(2\tau-1)(4\tau-1).$$
In view of (\ref{eq:a33}), 
\begin{equation}\label{eq:1a33}
a_{1,0,3,3}^\vee(1;\nu)=-1-18\mu-9\mu^2.
\end{equation}
In view of (\ref{eq:b33}), 
\begin{equation}\label{eq:1b33}
a_{1,0,3,3}^\wedge(1;\nu)=7+20\mu+\mu^2.
\end{equation}
Hence,
\begin{equation}\label{eq:as33}
a_{1,0,3,3}^\ast(1;\nu)=\tau^5-9\tau^4+20\tau^3-18\tau^2+7\tau-1=
\end{equation}
$$(\tau-1)(\tau^4-8\tau^3+12\tau^2-6\tau+1)=$$
$$(\tau-1)^2(\tau^3-7\tau^2+5\tau-1).$$
In view of (\ref{eq:a43}), (\ref{eq:15.1.b}), (\ref{eq:15.1.c}),
\begin{equation}\label{eq:a1v43}
a_{1,0,4,3}^{\vee1\vee}(1;\nu)=a_{1,0,4,3}^\vee(1;\nu)=
-1-31\mu-48\mu^2-4\mu^3,
\end{equation}
In view of (\ref{eq:b43}), (\ref{eq:15.1.b}), (\ref{eq:15.1.c}),
\begin{equation}\label{eq:b1w43}
a_{\alpha,0,4,3}^{\wedge1\wedge}(1;\nu)=a_{\alpha,0,4,3}^{\wedge}(1;\nu)=
9+53\mu+22\mu^2,
\end{equation}
In view of (\ref{eq:a1v43}), (\ref{eq:b1w43}), (\ref{eq:15c})
 and (\ref{eq:15.2.e}),
\begin{equation}\label{eq:a1s43}
a_{1,0,4,3}^{\ast1\ast}(1;\nu)=
\end{equation}
$$-4\tau^6+22\tau^5-48\tau^4+53\tau^3-31\tau^2+9\tau-1=$$
$$-(\tau-1)(4\tau^5-18\tau^4+30\tau^3-23\tau^2+8\tau-1)=$$
$$-(\tau-1)^2(4\tau^4-14\tau^3+16\tau^2-7\tau+1)=$$
$$
-(\tau-1)^3(4\tau^3-10\tau^2+6\tau-1)=-(\tau-1)^3(2\tau-1)(2\tau^2-4\tau+1).
$$
In view of (\ref{eq:a12}), 
\begin{equation}\label{eq:1a12}
a_{1,0,1,2}^\vee(1;\nu)=-20\mu-24\mu^2,
\end{equation}
In view of (\ref{eq:b12}), 
\begin{equation}\label{eq:1b12}
a_{1,0,1,2}^\wedge(1;\nu)=4+36\mu,
\end{equation}
In view of (\ref{eq:a22}), 
\begin{equation}\label{eq:1a22}
a_{1,0,2,2}^\vee(1;\nu)=-14\mu-41\mu^2.
\end{equation}
In view of (\ref{eq:b22}), 
\begin{equation}\label{eq:1b22}
a_{1,0,2,2}^\wedge(1;\nu)=2+36\mu+17\mu^2.
\end{equation}
In view of (\ref{eq:a32}), 
\begin{equation}\label{eq:1a32}
a_{1,0,3,2}^\vee(1;\nu)=-14\mu-58\mu^2-12\mu^3.
\end{equation}
In view of (\ref{eq:b32}), 
\begin{equation}\label{eq:1b32}
a_{1,0,3,2}^\wedge(1;\nu)=2+40\mu+42\mu^2.
\end{equation}
In view of (\ref{eq:a42}). 
\begin{equation}\label{eq:1a42}
a_{1,0,4,2}^\vee(1;\nu)=-15\mu-79\mu^2-38\mu^3.
\end{equation}
In view of (\ref{eq:b42}), 
\begin{equation}\label{eq:1b42}
a_{1,0,4,2}^\wedge(1;\nu)=2+47\mu+75\mu^2+8\mu^3.
\end{equation}
In view of (\ref{eq:a14}), (\ref{eq:b14}), (\ref{eq:a24}),
 (\ref{eq:b24}), (\ref{eq:a34}) and (\ref{eq:b34}),  
\begin{equation}\label{eq:1abk4}
a_{1,0,k,4}^\vee(1;\nu)=a_{1,0,k,4}^\wedge(1;\nu)=0.
\end{equation}
for $k=1,\,2,\,3.$
In view of (\ref{eq:a44}),
\begin{equation}\label{eq:1a44}
a_{1,0,4,4}^\vee(1;\nu)=-1-10\mu-5\mu^2.
\end{equation}
In view of (\ref{eq:b44}),
\begin{equation}\label{eq:1b44}
a_{1,0,4,4}^\wedge(1;\nu)=5+10\mu+\mu^2.
\end{equation}
In view of (\ref{eq:15.1.b}), (\ref{eq:15.1.c}), (\ref{eq:15.1.e})
and (\ref{eq:93bd}) with $\alpha=1,$
\begin{equation}\label{eq:av234,2}
a_{1,0,k,2+2j}^{\ast\kappa\ast}(1;\nu)=
-r_{1,0,2+2j}(\nu)a_{1,0,k,2\kappa+1,}^\ast(1;\nu)+
\end{equation}
$$r_{1,0,2\kappa+1}(\nu)a_{1,0,k,2+2j}^\ast(1;\nu)=
(-2)^\kappa\mu^{2-\kappa}a_{1,0,k,2+2j}^\ast(1;\nu),$$
where $j,\,\kappa=0,\,1,\,k=3-2\kappa,\,2,\,4.$
Therefore, in view of (\ref{eq:1a12}) -- (\ref{eq:1b44}),
\begin{equation}\label{eq:1as32}
a_{1,0,3,2}^{\ast0\ast}(1;\nu)=
\end{equation}
$$-\tau^5(12\tau^5-42\tau^4+58\tau^3-40\tau^2+14\tau-2)=$$
$$-\tau^5(\tau-1)(12\tau^4-30\tau^3+28\tau^2-12\tau+2)=$$
$$-\tau^5(\tau-1)^2(12\tau^3-18\tau^2+10\tau-2)=$$
$$-\tau^5(\tau-1)^2(2\tau-1)(6\tau^2-6\tau+2)=$$
$$-2\tau^5(\tau-1)^2(2\tau-1)(\tau^3-(\tau-1)^3),$$
\begin{equation}\label{eq:1a1s12}
a_{1,0,1,2}^{\ast1\ast}(1;\nu)=-2\tau^2(-24\tau^4+36\tau^3-20\tau^2+4\tau)=
\end{equation}
$$8\tau^3(6\tau^3-9\tau^2+5\tau-1)=8\tau^3(2\tau-1)(3\tau^2-3\tau+1)=$$
$$8\tau^3(2\tau-1)(\tau^3-(\tau-1)^3),$$
\begin{equation}\label{eq:1as22}
a_{1,0,2,2}^{\ast\kappa\ast}(1;\nu)=
(-2)^\kappa\tau^{4-2\kappa}(17\tau^5-41\tau^4+36\tau^3-14\tau^2+2\tau)=
\end{equation}
$$(-2)^\kappa\tau^{5-2\kappa}(\tau-1)(17\tau^3-24\tau^2+12\tau-2)=$$
$$(-2)^\kappa\tau^{5-2\kappa}(\tau-1)(\tau^3+2(2\tau-1)^3),$$
\begin{equation}\label{eq:1as42}
a_{1,0,4,2}^{\ast\kappa\ast}(1;\nu)=
\end{equation}
$$(-2)^\kappa\tau^{5-2\kappa}
(8\tau^6-38\tau^5+75\tau^4-79\tau^3+47\tau^2-15\tau+2)=$$
$$(-2)^\kappa\tau^{5-2\kappa}
(\tau-1)(8\tau^5-30\tau^4+45\tau^3-34\tau^2+13\tau-2)=$$
$$(-2)^\kappa\tau^{5-2\kappa}(\tau-1)^2(8\tau^4-22\tau^3+23\tau^2-11\tau+2)=$$
$$(-2)^\kappa\tau^{5-2\kappa}(\tau-1)^3(8\tau^3-14\tau^2+9\tau-2)=$$
$$(-2)^\kappa\tau^{5-2\kappa}(\tau-1)^3(2\tau-1)(4\tau^2-5tau+2).$$
\begin{equation}\label{eq:vvwwk4}
a_{1,0,k,4}^{\ast\kappa\ast}(1;\nu)=0\,\,\text{for}\,\,k=1,\,2,\,3,
\end{equation}
\begin{equation}\label{eq:1as44}
a_{1,0,4,4}^{\ast\kappa\ast}(1;\nu)=
\end{equation}
$$(-2)^\kappa\tau^{4-2\kappa}(\tau^5-5\tau^4+10\tau^3-10\tau^2+5\tau-1)=$$
$$(-2)^\kappa\tau^{4-2\kappa}(\tau-1)^5,$$
where $\kappa=0,\,1.$
Let
\begin{equation}\label{eq:v2}
\theta(k,\kappa)=\det\left(\matrix
r_{1,0,1+2\kappa}&r_{1,0,3-2\kappa}\\
a_{1,0,k,2\kappa+1}^\ast(1;\nu)&a_{1,0,k,3-2\kappa}^\ast(1;\nu),
\endmatrix\right),
\end{equation}
where $\kappa=0,\,1,\,k=1,\,2,\,3,\,4.$
Then
\begin{equation}\label{eq:v3}
 \theta(k,\kappa)=-\theta(k,1-\kappa),\,\,\text{for}
\,\,\kappa=0,\,1,\,k=1,\,2,\,3,\,4;
\end{equation}
 moreover,
in view of (\ref{eq:93bd}) with $\alpha=1$, (\ref{eq:15.1.b}),
 (\ref{eq:15.1.c}) and (\ref{eq:15.1.e}),
\begin{equation}\label{eq:av234,3}
a_{1,0,k,3-2\kappa}^{\ast\kappa\ast}(1;\nu)=\theta(k,\kappa)=
-r_{1,0,3-2\kappa)}(\nu)
a_{1,0,k,2\kappa+1}^\ast(1;\nu)+
\end{equation}
$$r_{1,0,1+2\kappa}(\nu)a_{1,0,k,3-2\kappa}^\ast(1;\nu)=
-(-2)^{1-\kappa}\mu^{1+\kappa}a_{1,0,k,2\kappa+1}^\ast(1;\nu)+$$
$$(-2)^\kappa\mu^{2-\kappa} a_{1,0,k,3-2\kappa}^\ast(1;\nu),$$
where $\kappa=0,\,1,\,k=3-2\kappa,\,2,\,4.$
If we replace $\kappa$ by $\kappa^\prime=1-\kappa$ in (\ref{eq:av234,3}),
then we obtain the equality 
\begin{equation}\label{eq:avh234,3}
a_{1,0,k,3-2\kappa^\prime}^{\ast\kappa^\prime\ast}(1;\nu)=
\theta(k,\kappa^\prime)
\end{equation}
where $\kappa^\prime=1-\kappa=0,\,1,\,k=3-2\kappa^\prime=2\kappa+1,\,2,\,4.$
Hence,
\begin{equation}\label{eq:v4}
a_{1,0,k,3-2\kappa}^{\ast\kappa\ast}(1;\nu)=
-a_{1,0,k,2\kappa+1}^{\ast1-\kappa\ast}(1;\nu),
\end{equation}
where $\kappa=0,\,1,\,k=2,\,4.$
Therefore, in view of (\ref{eq:as11}) -- (\ref{eq:a1s43}),
\begin{equation}\label{eq:1as11}
a_{1,0,1,1}^{\ast1\ast}(1;\nu)=-2\tau^4(17\tau^3-21\tau^2+10\tau-2)-
\end{equation}
$$-2\tau^4(-6\tau^2+5\tau-1)=-2\tau^4(17\tau^3-27\tau^2+15\tau-3)=$$
$$-2\tau^4((\tau-1)^3+2(2\tau-1)^3).$$
\begin{equation}\label{eq:1as33}
a_{1,0,3,3}^{\ast0\ast}(1;\nu)=
\end{equation}
$$\tau^4(\tau-1)^2(\tau^3-7\tau^2+5\tau-1)+$$
$$2\tau^2(\tau^2(\tau-1)^2(8\tau^3-10\tau^2+5\tau-1))=$$
$$\tau^4(\tau-1)^2(17\tau^3-27\tau^2+15\tau-3)=$$
$$\tau^4(\tau-1)^2((\tau-1)^3+2(2\tau-1)^3)=$$
$$-\frac12(\tau-1)^2a_{1,0,1,1}^{\ast1\ast}(1;\nu).$$
\begin{equation}\label{eq:1as23}
a_{1,0,2,3}^{\ast0\ast}(1;\nu)=-a_{1,0,2,1}^{\ast1\ast}(1;\nu)=
\end{equation}
$$\tau^4((\tau-1)(2\tau-1)(4\tau-1))+$$
$$2\tau^2(-\tau^2(\tau-1)(2\tau-1)(6\tau^2-4\tau+1)=$$
$$-\tau^4(\tau-1)(2\tau-1)(12\tau^2-12\tau+3)=$$
$$-3\tau^4(\tau-1)(2\tau-1)^3.$$
\begin{equation}\label{eq:1as43}
a_{1,0,4,3}^{\ast0\ast}(1;\nu)=-a_{1,0,4,1}^{\ast1\ast}(1;\nu)=
\end{equation}
$$\tau^4(-(\tau-1)^3(2\tau-1)(2\tau^2-4\tau+1))+$$
$$2\tau^2(-\tau^2(\tau-1)^3(2\tau-1)(2\tau^2-2\tau+1))=$$
$$-\tau^4(\tau-1)^3(2\tau-1)(6\tau^2-8\tau+3).$$
{\begin{center}\large\bf
\S 5. Further properties of the matrix $A^{\ast\kappa\ast}_{1,0}(z,\nu).$
\end{center}}
In view of (\ref{eq:15.2.a11}), (\ref{eq:1as44}) 
\begin{equation}\label{eq:ast44,11,33}
\lim\limits_{\nu\to+\infty}(-2)^{-\kappa} \nu^{-9+2\kappa}
a_{1,0,4,4}^{\vee\kappa\vee}(1;\nu)=
\end{equation}
$$\lim\limits_{\nu\to+\infty}(-2)^{-\kappa} \nu^{-9+2\kappa}
a_{1,0,2\kappa+1,2\kappa+1}^{\vee\kappa\vee}(1;\nu)=1,$$
In view of (\ref{eq:15.2.vw1,234})
\begin{equation}\label{eq:ast1k,3k}
\lim\limits_{\nu\to+\infty}
(-2)^\kappa\nu^la_{1,0,2\kappa+1,k}^{\ast\kappa\ast}(1;\nu)=0,
\end{equation}
where $\kappa=0,\,1,\,k=3-2\kappa,\,2,\,4,\,l\in{\mathbb Z}.$
In view of (\ref{eq:as21}) and (\ref{eq:a1s23}),
\begin{equation}\label{eq:ast21,23}
\lim\limits_{\nu\to+\infty}
(-2)^{-\kappa} \nu^{-6}a_{1,0,2,2\kappa+1}^{\ast\kappa\ast}(1;\nu)=-12(1-\kappa,
\end{equation}
where $\kappa=0,\,1.$
In view of (\ref{eq:1as22}), (\ref{eq:1as11}) (\ref{eq:1as22})
 and (\ref{eq:1as33}),
\begin{equation}\label{eq:ast22,33,11}
\lim\limits_{\nu\to+\infty}(-2)^{-\kappa} \nu^{-9+2\kappa}
a_{1,0,2,2}^{\vee\kappa\vee}(1;\nu)=
\end{equation}
$$\lim\limits_{\nu\to+\infty}(-2)^{-\kappa} \nu^{-9+2\kappa}
a_{1,0,3-2\kappa,3-2\kappa}^{\vee\kappa\vee}(1;\nu)=17,$$
where $\kappa=0,\,1.$
In view of (\ref{eq:1as23}),
\begin{equation}\label{eq:ast23,21}
\lim\limits_{\nu\to+\infty}
(-2)^{-\kappa}\nu^{-8}a_{1,0,2,3-2\kappa}^{\ast\kappa\ast}(1;\nu)=
-2^{-\kappa}24,
\end{equation}
where $\kappa=0,\,1.$
In view of (\ref{eq:vvwwk4}),
\begin{equation}\label{eq:15.3.astk4}
\lim\limits_{\nu\to+\infty}(-2)^{-\kappa} \nu^la_{1,0,k,4}^{\ast\ast}(1;\nu)=0,
\end{equation}
where $\kappa=0,\,1,\,k=1\,2,\,3,\,l\in{\mathbb Z}.$
In view of (\ref{eq:as31}), (\ref{eq:a1s13}),
\begin{equation}\label{eq:ast31,13}
\lim\limits_{\nu\to+\infty}(-2)^{-\kappa}\nu^{-7+4\kappa}
a_{1,0,3-2\kappa,1+2\kappa}^{\ast\kappa\ast}(1;\nu)=8(1-\kappa)
\end{equation}
In view of (\ref{eq:1as32}) and (\ref{eq:1a1s12})
\begin{equation}\label{eq:ast32,12}
\lim\limits_{\nu\to+\infty} (-2)^{-\kappa}\nu^{-10+4\kappa}
a_{1,0,3-2\kappa,2}^{\ast\ast}(1;\nu)=
-2^\kappa12, 
\end{equation}
where $\kappa=0,\,1.$
In view of (\ref{eq:as41}) and (\ref{eq:a1s43})
\begin{equation}\label{eq:15.3.ast41,43}
\lim\limits_{\nu\to+\infty} (-2)^{-\kappa}\nu^{-8+2\kappa}
a_{1,0,4,2\kappa+1}^{\ast\kappa\ast}(1;\nu)=-(-2)^{-\kappa}4, 
\end{equation}
where $\kappa=0,\,1.$
In view of (\ref{eq:1as42}),
\begin{equation}\label{eq:15.3.ast42}
\lim\limits_{\nu\to+\infty}
(-2)^{-\kappa}\nu^{-11+2\kappa}a_{1,0,4,2}^{\ast\kappa\ast}(1;\nu)=
8, 
\end{equation}
where $\kappa=0,\,1.$
In view of (\ref{eq:1as43}),
\begin{equation}\label{eq:ast43,41}
\lim\limits_{\nu\to+\infty} (-2)^{-\kappa}\nu^{-10}
a_{1,0,4,3-2\kappa}^{\ast\kappa\ast}(1;\nu)=-2^{-\kappa}12, 
\end{equation}
where $\kappa=0,\,1.$
 
Let $\lambda$ be a variable. We denote by $T_{n,\lambda}$ the
diagonal $n\times n$-matrix, $i$-th diagonal element of which
is equal to $\lambda^{i-1}$ for $i=1,\,...,\,n.$
We denote by $T^{(i)}_{\nu,\lambda}$ $n\times n$-diagnoal matrix such that
 its $i$-th  diagonal element is equal to $\lambda$
 and and its other diagonal elements
are equal to $1.$
It follows from (\ref{eq:ast44,11,33}) -- (\ref{eq:ast43,41})
that
\begin{equation}\label{eq:15.3.h}
\lim\limits_{\nu\to+\infty}
(-2)^{-\kappa}\nu^{-9+2\kappa}
((T^{(2\kappa+1)}_{4,\nu})^{4-2\kappa}T_{4,\nu})^{-1}
A^{\ast\kappa\ast}_{1,0}(1;\nu)
(T^{(2\kappa+1)}_{4,\nu})^{4-2\kappa}T_{4,\nu})
=A^{(\kappa)},
\end{equation}
where $\kappa=0,1,$
$$A^{(0)}=\left(\matrix
  1&  0&  0&0\\
-12& 17&-24&0\\
  8&-12& 17&0\\
  4&  8&-12&1
\endmatrix\right),
A^{(1)}=\left(\matrix
 17&-24&0&0\\
-12& 17&0&0\\
  0&  0&1&0\\
 -6&  8&2&1
\endmatrix\right).$$
Let me to make some comments to the equalities (\ref{eq:15.3.h}.
Accordng to (\ref{eq:ast1k,3k}) and (\ref{eq:15.3.h} 
\begin{equation}\label{eq:15.3.h0}
\lim\limits_{\nu\to+\infty}
\nu^{-9}((T^{(1)}_{4,\nu})^{4-l}T_{4,\nu})^{-1}
A^{\ast\kappa\ast}_{1,0}(1;\nu)T^{(1)}_{4,\nu})^{4-l}T_{4,\nu})
\end{equation}
$$=\left(\matrix
  1&  0&  0&0\\
  0& 17&-24&0\\
  0&-12& 17&0\\
  04&  8&-12&1
\endmatrix\right),$$
for any $l\in{\mathbb N},$ and, in particular, for $l=4.$
 Let $\Sigma_\kappa$ be the matrix, which appears
after replacement $2\kappa+1$-th row in $E_4$ by the row
$\left(\matrix1&0&-2&0\endmatrix\right).$
 Then, in view of (\ref{eq:15.1.a})
with $\alpha=1$ and (\ref{eq:15.1.b},  
$$T_{4,\nu}^{-1}C^{(\kappa)}(\nu)T_{4,\nu}=
(T^{2\kappa+1}_{4,\nu})^{4-2\kappa}\Sigma_\kappa.$$
Further we have
$$(-2)^{-\kappa}\nu^{-9+2}\times$$
$$((T^{(2\kappa+1}_{4,\nu})^{4-2\kappa}T_{4,\nu})^{-1}
A^{\ast\kappa\ast}_{1,0}(1,\nu)T_{4,\nu}
T^{(2\kappa+1}_{4,\nu})^{4-2\kappa}=$$
$$(T^{(2\kappa+1}_{4,\nu})^{2\kappa-4}\times$$
$$(T_{4,\nu}^{-1}C^{(\kappa)}(\nu-1)T_{4,\nu})\times$$
$$T_{4,\nu}^{-1}A^\ast_{1,0}(1,\nu)T_{4,\nu}\times$$
$$(T_{4,\nu}^{-1}(C^{(\kappa)})^{-1}T_{4,\nu})=$$
$$
\Sigma_\kappa T_{4,\nu}^{-1}A^\ast_{1,0}(1,\nu)T_{4,\nu}(\Sigma_\kappa)^{-1}.
$$
According to the results of section 11.1 in [\ref{r:caj13}],
$$\lim\limits_{\nu\to\infty}T_{4,\nu}^{-1}A^\ast_{1,0}(1,\nu)T_{4,\nu}
=A^{\sim}=
$$
$$\left(\matrix
 17&-24&0&0\\
-12& 17&0&0\\
  8&-12&1&0\\
 -4&  8&-4&1
\endmatrix\right).$$
Further we have
$$\Sigma_0A^{\sim}(\Sigma_0)^{-1}=\left(\matrix
  1&  0&-2&0\\
-12& 17& 0&0\\
  8&-12& 1&0\\
 -4&  8&-4&1
\endmatrix\right)\left(\matrix
  1&0&2&0\\
  0&1&0&0\\
  1&0&1&0\\
  0&0&0&1
\endmatrix\right)=$$
$$\left(\matrix
  1&  0&  0&0\\
-12& 17&-24&0\\
  8&-12& 17&0\\
 -4&  8& -12&1
\endmatrix\right)=A^{(0)},
$$
$$\Sigma_1A^{\sim}(\Sigma_1)^{-1}=
\left(\matrix
 17&-24&0&0\\
-12& 17&0&0\\
  1&  0&-2&0\\
 -4&  8&-4&1
\endmatrix\right)\left(\matrix
      1&0&         0&0\\
      0&1&         0&0\\
\frac12&0&-\frac12&0\\
      0&0&          0&1
\endmatrix\right)=$$
$$\left(\matrix
 17&-24&0&0\\
-12& 17&0&0\\
  0&  0&1&0\\
 -6&  8&2&1
\endmatrix\right)=A^{(1)}.$$
So, I can consider the equalities (\ref{eq:15.3.h})
as some test of previous calculations.

{\begin{center}\large\bf
\S 6. Further properties of the functions considered in
 the Introduction.
\end{center}}
Let $\alpha\ge1.$ Then the function $t^r(R(\alpha,t,\nu))^2$ (see (\ref{eq:d}))
 is regular at $t=\infty$ for $r=0,\,1,\,2,$ is regular at $t=\infty$
 for $r=3,\,\alpha\ge2$ and has a pole of first order at $t=\infty$
 for $r=3,\,\alpha=1.$ So, in the case $r=0,\,1,\,2,\,\alpha\ge1$ we 
have the equalities
\begin{equation}\label{eq:15.4.a.0a}
\Res(t^r(R(\alpha,t,\nu))^2,t=\infty)=0,\,
\end{equation}
\begin{equation}\label{eq:15.4.a.0b}
\lim\limits_{t\to \infty} t^r(R(\alpha,t,\nu))^2=0,
\end{equation}
and in the case $r=3,\,\alpha=1$ we 
have the equalities
\begin{equation}\label{eq:15.4.a.1a}
\Res(t^3(R(1,t,\nu))^2,t=\infty)=-1,\,
\end{equation}
\begin{equation}\label{eq:15.4.a.1b}
\lim\limits_{t\to \infty} t^3(R(1,t,\nu))^2=0.
\end{equation}
 In view of (\ref{eq:h}),
\begin{equation}\label{eq:15.4.hr}
\delta^rf_{\alpha,0,2}^{\ast}(z, \nu)=
\sum\limits_{t=1}^{+\infty}z^{-t}
\frac{(\nu+\alpha)!^2}{\nu!^2}(-t)^r(R(\alpha,t,\nu))^2,
\end{equation}
where we consider $r=0,\,1,\,2,\,3.$
Expanding  $\frac{(\nu+\alpha)!^2}{(\nu!)^2}(-t)^r(R(\alpha,t,\nu))^2$
 into partial fractions relatively $t$, we obtain
\begin{equation}\label{eq:6bg}
\frac{(\nu+\alpha)!^2}{(\nu!)^2}(-t)^r(R(\alpha,t;\nu))^2=
\end{equation}
$$\sum\limits_{i=1}^2
\left(\sum\limits_{k=0}^{\nu+\alpha}
\beta_{\alpha,0,i,k,\nu}^{(r)}(t + k)^{-i}\right),$$
where $\nu\in{\mathbb N}_0,\,r=0,\,1,\,2,\,3,\,\alpha\in{\mathbb N},$
\begin{equation}\label{eq:6bh}
\beta_{\alpha,0,2-j,k,\nu}^{(r)}=\frac{(\nu+\alpha)!^2}{(\nu!)^2}\times
\end{equation}
$$\frac1{j!}\lim\limits_{t\to-k}\left(\frac{\partial}{\partial t}\right)^j
((-t)^r(R(\alpha,t,\nu)(t+k))^2),$$ for $j=0,\,1.$
In view of (\ref{eq:15.4.a.0a}), (\ref{eq:6bg})  and (\ref{eq:15.4.a.1a}),
\begin{equation}\label{eq:6bh.0}
\sum\limits_{k=0}^{\nu+\alpha}\beta_{\alpha,0,1,k,\nu}^{(r)}=0,
\end{equation}
 for $r=0,\,1,\,2$ and $\alpha\in{\mathbb N}.$
\begin{equation}\label{eq:6bh.1}
\sum\limits_{k=0}^{\nu+1}\beta_{1,0,1,k,\nu}^{(3)}=-(\nu+1)^2.
\end{equation}
In view of (\ref{eq:6bg}),
\begin{equation}\label{eq:6bgf}
-\frac{(\nu+\alpha)!^2}{(\nu!)^2}\frac\partial
{\partial t}((-t)^r(R(\alpha,t;\nu))^2)=
\end{equation}
$$\sum\limits_{i=1}^2
\left(\sum\limits_{k=0}^{\nu+\alpha}
\beta_{\alpha,0,i,k,\nu}^{(r)}i(t + k)^{-i-1}\right),$$
where $\nu\in{\mathbb N}_0,\,r=0,\,1,\,2,\,3.$
 Let
\begin{equation}\label{eq:6bi0}
S_{i,k}(\alpha,\nu)=-\left(\sum\limits_{\kappa=k+1}^{\nu+k}1/\kappa^i\right)-
\end{equation}
$$\left(\sum\limits_{\kappa=1}^{\nu+\alpha-k}1/\kappa^i\right)+
\sum\limits_{\kappa=1}^k1/\kappa^i,$$
where 
$\nu\in{\mathbb N}_0,\,i\in{\mathbb N},\,k\in[0,\,\nu+\alpha\cap{\mathbb Z}.$
In particular,
\begin{equation}\label{eq:6bi1a0}
S_{1,0}(0,0)=0,\,S_{1,0}(0,1)=-2,\,S_{1,1}(0,1)=\frac12,
\end{equation}
\begin{equation}\label{eq:6bi1k01a1n0}
S_{1,0}(1,0)=-1,\,S_{1,1}(1,0)=1,
\end{equation}
\begin{equation}\label{eq:6bi1k012a1n1}
S_{1,0}(1,1)=-\frac52,\,S_{1,1}(1,1)=-\frac12,\,
S_{1,2}(1,1)=\frac76,
\end{equation}
\begin{equation}\label{eq:6bi1k0a1n2}
S_{1,0}(1,2)=-(1+1/2)-(1+1/2+1/3)=-\frac{10}3
\end{equation}
\begin{equation}\label{eq:6bi1k1a1n2}
S_{1,1}(1,2)=-(1/2+1/3)-(1+1/2)+1=-\frac43
\end{equation}
\begin{equation}\label{eq:6bi1k2a1n2}
S_{1,2}(1,2)=-(1/3+1/4)-1+(1+1/2)=-\frac1{12}
\end{equation}
\begin{equation}\label{eq:6bi1k3a1n2}
S_{1,3}(1,2)=-(1/4+1/5)+(1+1/2+1/3)=\frac{83}{60}
\end{equation}
In view of (\ref{eq:6bh}), (\ref{eq:d}) and (\ref{eq:6bi0})
\begin{equation}\label{eq:6bi2}
\beta_{\alpha,0,2,k,\nu}^{(0)}=
\end{equation}
$$ \left(\frac{(\nu+\alpha)!}{\nu!}\times\frac
{(\nu+k)!}{k!}\times\frac1{(\nu+\alpha-k)!}\right)^2=
\binom{\nu+\alpha}k^2\binom{\nu+k}k^2,$$
\begin{equation}\label{eq:6bi1}
\beta_{\alpha,0,1,k,\nu}^{(0)}=2\beta_{\alpha,0,2,k,\nu}^{(0)}
S_{1,k}(\alpha,\nu),
\end{equation}
where 
$\nu\in{\mathbb N}_0,\,i\in{\mathbb N},\,k\in[0,\,\nu+\alpha\cap{\mathbb Z}.$
In particular,
\begin{equation}\label{eq:r0a0i2k01n01}
\beta_{0,0,2,0,0}^{(0)}=\beta_{0,0,2,0,1}^{(0)}=1,\,\
\beta_{0,0,2,1,1}^{(0)}=4,
\end{equation}
\begin{equation}\label{eq:r0a1i2k01n0}
\beta_{1,0,2,0,0}^{(0)}=\beta_{1,0,2,1,0}^{(0)}=1,
\end{equation}
\begin{equation}\label{eq:r0a1i2k012n1}
\beta_{1,0,2,0,1}^{(0)}=1,\,\beta_{1,0,2,1,1}^{(0)}=16,\,
\beta_{1,0,2,2,1}^{(0)}=9,
\end{equation}
\begin{equation}\label{eq:r0a1i2k01n2}
\beta_{1,0,2,0,2}^{(0)}=1,\,\beta_{1,0,2,1,2}^{(0)}=81,
\end{equation}
\begin{equation}\label{eq:r0a1i2k23n2}
\beta_{1,0,2,2,2}^{(0)}=324,\,\beta_{1,0,2,3,2}^{(0)}=100.
\end{equation}
In view of (\ref{eq:6bi1}), (\ref{eq:r0a0i2k01n01}) -- (\ref{eq:r0a1i2k23n2})
 and (\ref{eq:6bi1a0}) -- (\ref{eq:6bi1k3a1n2}),
\begin{equation}\label{eq:r0a0i1k0n0}
\beta_{0,0,1,0,0}^{(0)}=2\beta_{0,0,2,0,0}^{(0)}S_{1,0}(0,0)=0,
\end{equation}
\begin{equation}\label{eq:r0a0i1k0n1}
\beta_{0,0,1,0,1}^{(0)}=2\beta_{0,0,2,0,1}^{(0)}S_{1,0}(0,1)=
2\times1\times(-2)=-4,
\end{equation}
\begin{equation}\label{eq:r0a0i1k1n1}
\beta_{0,0,1,1,1}^{(0)}=
2\beta_{0,0,2,1,1}^{(0)}S_{1,1}(0,1)=2\times4\times\frac12=4,
\end{equation}
\begin{equation}\label{eq:r0a1i1k0n0}
\beta_{1,0,1,0,0}^{(0)}=2\beta_{1,0,2,0,0}^{(0)}S_{1,0}(1,0)=2\times1\times(-1)=-2,
\end{equation}
\begin{equation}\label{eq:r0a1i1k1n0}
\beta_{1,0,1,1,0}^{(0)}=
2\beta_{1,0,2,1,0}^{(0)}S_{1,1}(1,0)=2\times1\times1=2,
\end{equation}
\begin{equation}\label{eq:r0a1i1k0n1}
\beta_{1,0,1,0,1}^{(0)}=
2\beta_{1,0,2,0,1}^{(0)}S_{1,0}(1,1)=2\times1\times\frac{-5}2=-5,
\end{equation}
\begin{equation}\label{eq:r0a1i1k1n1}
\beta_{1,0,1,1,1}^{(0)}=2\beta_{1,0,2,1,1}^{(0)}S_{1,1}(1,1)=
2\times16\times\frac{-1}2=-16,
\end{equation}
\begin{equation}\label{eq:r0a1i1k2n1}
\beta_{1,0,1,2,1}^{(0)}=2\beta_{1,0,2,2,1}^{(0)}S_{1,2}(1,1)=
2\times9\times\frac76=21,
\end{equation}
\begin{equation}\label{eq:r0a1i1k0n2}
\beta_{1,0,1,0,2}^{(0)}=2\beta_{1,0,2,0,2}^{(0)}S_{1,0}(1,2)=
2\times1\times\frac{-10}3=\frac{-20}3,
\end{equation}
\begin{equation}\label{eq:r0a1i1k1n2}
\beta_{1,0,1,1,2}^{(0)}=2\beta_{1,0,2,1,2}^{(0)}S_{1,1}(1,2)=
2\times81\times\frac{-4}3=-216,
\end{equation}
\begin{equation}\label{eq:r0a1i1k2n2}
\beta_{1,0,1,2,2}^{(0)}=2\beta_{1,0,2,2,2}^{(0)}S_{1,2}(1,2)=
2\times324\times\frac{-1}{12}=-54,
\end{equation}
\begin{equation}\label{eq:r0a1i1k3n2}
\beta_{1,0,1,3,2}^{(0)}=2\beta_{1,0,2,3,2}^{(0)}S_{1,3}(1,2)=
2\times100\times\frac{83}{60}=\frac{830}3.
\end{equation}
We put in (\ref{eq:6bg}) $r=0,$ and multiply both sides of obtained
equality by $(-t)^r$ for $r=0,\,1,\,2,\,3.$ Then we see that
\begin{equation}\label{eq:6bg.1}
-t\frac{(\nu+\alpha)!^2}{(\nu!)^2}(R(\alpha,t;\nu))^2=
\end{equation}
$$\sum\limits_{i=1}^2
\left(\sum\limits_{k=0}^{\nu+\alpha}
\frac{\beta_{\alpha,0,i,k,\nu}^{(0)}(-t-k+k)}{(t + k)^i}\right)=
\left(\sum\limits_{k=0}^{\nu+\alpha}
\frac{k\beta_{\alpha,0,2,k,\nu}^{(0)}}{(t + k)^2}\right)+$$
$$\left(\sum\limits_{k=0}^{\nu+\alpha}
\frac{k\beta_{\alpha,0,1,k,\nu}^{(0)}-\beta_{\alpha,0,2,k,\nu}^{(0)}}
{t + k}\right)-
\sum\limits_{k=0}^{\nu+\alpha}
\beta_{\alpha,0,1,k,\nu}^{(0)},$$
\begin{equation}\label{eq:6bg.2}
(-t)^2\frac{(\nu+\alpha)!^2}{(\nu!)^2}(R(\alpha,t;\nu))^2=
\end{equation}
$$\sum\limits_{i=1}^2
\left(\sum\limits_{k=0}^{\nu+\alpha}
\frac{\beta_{\alpha,0,i,k,\nu}^{(0)}(t+k-k)^2}{(t + k)^i}\right)=
\left(\sum\limits_{k=0}^{\nu+\alpha}
\frac{k^2\beta_{\alpha,0,2,k,\nu}^{(0)}}{(t + k)^2}\right)+$$
$$\left(\sum\limits_{k=0}^{\nu+\alpha}\frac{
k^2\beta_{\alpha,0,1,k,\nu}^{(0)}-2k\beta_{\alpha,0,2,k,\nu}^{(0)}}
{t+k}\right)+
\sum\limits_{k=0}^{\nu+\alpha}(
\beta_{\alpha,0,2,k,\nu}^{(0)}+(t-k)\beta_{\alpha,0,1,k,\nu}^{(0)}
),$$
\begin{equation}\label{eq:6bg.3}
(-t)^3\frac{(\nu+\alpha)!^2}{(\nu!)^2}(R(\alpha,t;\nu))^2=
\end{equation}
$$\sum\limits_{i=1}^2
\left(\sum\limits_{k=0}^{\nu+\alpha}
\frac{\beta_{\alpha,0,i,k,\nu}^{(0)}(-t-k+k)^3}{(t + k)^i}\right)=
\left(\sum\limits_{k=0}^{\nu+\alpha}\frac
{k^3\beta_{\alpha,0,2,k,\nu}^{(0)}}{(t + k)^2}\right)+$$
$$\left(\sum\limits_{k=0}^{\nu+\alpha}\frac
{k^3\beta_{\alpha,0,1,k,\nu}^{(0)}-3k^2\beta_{\alpha,0,2,k,\nu}^{(0)}}
{t + k}\right)-
\left(\sum\limits_{k=0}^{\nu+\alpha}
(t-2k)(\beta_{\alpha,0,2,k,\nu}^{(0)}\right)-$$
$$\sum\limits_{k=0}^{\nu+\alpha}
(t^2-kt+k^2)\beta_{\alpha,0,1,k,\nu}^{(0)}.$$
The equality (\ref{eq:6bh.0}) with $r=0$ again follows
 from (\ref{eq:15.4.a.0b}) with $r=1$ and (\ref{eq:6bg.1}); moreover,
in view of (\ref{eq:6bg}) with $r=1,$ and (\ref{eq:6bg.1}), 
\begin{equation}\label{eq:6bg.1.2}
\beta_{\alpha,0,2,k,\nu}^{(1)}=k\beta_{\alpha,0,2,k,\nu}^{(0)},
\end{equation}
\begin{equation}\label{eq:6bg.1.1}
\beta_{\alpha,0,1,k,\nu}^{(1)}=k\beta_{\alpha,0,1,k,\nu}^{(0)}-
\beta_{\alpha,0,2,k,\nu}^{(0)}
\end{equation}
for $k=0,\,...,\,\nu+\alpha,\,\alpha\in{\mathbb N}.$
The equality (\ref{eq:6bh.0}) with $r=1$ again follows
 from (\ref{eq:15.4.a.0b}) with $r=2,$ (\ref{eq:6bh.0}) with $r=0,$ 
(\ref{eq:6bg.2}) and (\ref{eq:6bg.1.1}); moreover,
in view of (\ref{eq:6bg}) with $r=2,$ and (\ref{eq:6bg.2}),
\begin{equation}\label{eq:6bg.2.2}
\beta_{\alpha,0,2,k,\nu}^{(2)}=k^2\beta_{\alpha,0,2,k,\nu}^{(0)},
\end{equation}
\begin{equation}\label{eq:6bg.2.1}
\beta_{\alpha,0,1,k,\nu}^{(2)}=k^2\beta_{\alpha,0,1,k,\nu}^{(0)}-
2k\beta_{\alpha,0,2,k,\nu}^{(0)}
\end{equation}
for $\alpha\in{\mathbb N}$ and $k=0,\,...,\,\nu+\alpha.$
The equality (\ref{eq:6bh.0}) with $r=2$ again follows
 from  (\ref{eq:15.4.a.1b}), (\ref{eq:6bh.0}) with both $r\in\{0,\,1\},$ 
(\ref{eq:6bg.3}), (\ref{eq:6bg.1.1}) and from (\ref{eq:6bg.2.1}); moreover,
in view of (\ref{eq:6bg}) with $r=3,$ and (\ref{eq:6bg.3}),
\begin{equation}\label{eq:6bg.3.2}
\beta_{\alpha,0,2,k,\nu}^{(3)}=k^3\beta_{\alpha,0,2,k,\nu}^{(0)},
\end{equation}
\begin{equation}\label{eq:6bg.3.1}
\beta_{\alpha,0,1,k,\nu}^{(3)}=k^3\beta_{\alpha,0,1,k,\nu}^{(0)}-
3k^2\beta_{\alpha,0,2,k,\nu}^{(0)}
\end{equation}
for $\alpha\in{\mathbb N}$ and $k=0,\,...,\,\nu+\alpha.$
In view of (\ref{eq:h}) -- (\ref{eq:ac}),
\begin{equation}\label{eq:acrf}
(\delta^r)f_{\alpha,0,3}^{\ast}(z,\nu)=
(\log(z))(\delta)^r)f_{\alpha,0,2}^{\ast}(z,\nu)+
\end{equation}
$$r(\delta)^{r-1}f_{\alpha,0,2}^{\ast}(z,\nu)
+(\delta)^rf_{\alpha,0,4}^{\ast}(z,\nu)=
(\log(z))(\delta)^r)f_{\alpha,0,2}^{\ast}(z,\nu)+$$
$$\sum\limits_{t=1}^{+\infty} z^{-t}
\frac{(\nu+\alpha)!^2}{\nu!^2}
\left (r(-t)^{r-1}-(-t)^r\frac\partial{\partial t}\right)R^2(\alpha,t,\nu)=$$
$$(\log(z))(\delta)^r)f_{\alpha,0,2}^{\ast}(z,\nu)-
\sum\limits_{t=1}^{+\infty} z^{-t}
\frac{(\nu+\alpha)!^2}{\nu!^2}
\frac\partial{\partial t}((-t)^rR^2)(\alpha,t,\nu).$$ 
In view of (\ref{eq:6bg}), (\ref{eq:6bgf}), (\ref{eq:15.4.hr})
(\ref{eq:acrf}),
\begin{equation}\label{eq:15.4.hr1f}
\delta^rf_{\alpha,0,2+j}^{\ast}(z, \nu)-
j(\log(z))\delta^rf_{\alpha,0,2}^{\ast}(z, \nu)=
\end{equation}
$$\sum\limits_{i=1}^2
\left(\sum\limits_{t=1}^\infty
\left(\sum\limits_{k=0}^{\nu+\alpha}(1-j+ij)\beta_{\alpha,0,i,k,\nu}^{(r)}
z^kz^{-t-k}(t + k)^{-i-j}\right)\right)=
$$
$$\sum\limits_{i=1}^2\left(\sum\limits_{k=0}^{\nu+\alpha}
(1-j+ij)\beta_{\alpha,0,i,k,\nu}^{(r)}z^k
\left(\sum\limits_{t=1}^\infty
z^{-t-k} (t + k)^{-i-j}\right)\right)=
$$
$$\sum\limits_{i=1}^2\left(\sum\limits_{k=0}^\nu
(1-j+ij)\beta_{\alpha,0,i,k,\nu}^{(r)}z^k\left(L_{i+1}(1/z)-
\sum\limits_{\tau=1}^kz^{-\tau}i(\tau)^{-i-j}\right)\right)=
$$
$$\left(\sum\limits_{i=1}^2
(1-j+ij)\beta^{\ast(r)}_{\alpha,0,i}(z;\nu)L_{i+j}(1/z)\right)-
\beta^{(r)}_{\alpha,0,3+j}(z;\nu),$$
where $j=0,\,1,\,r=0,\,1,\,2,\,3,\,\vert z\vert>1,\,\alpha\in{\mathbb N},$
\begin{equation}\label{eq:6cd}
L_s(1/z)=\sum\limits_{n=1}^\infty 1/(z^nn^s),\,
\beta^{\ast(r)}_{\alpha,0,i}(z;\nu)=
\sum\limits_{k=0}^{\nu+\alpha}\beta_{\alpha,0,i,k,\nu}^{(r)}z^k,
\end{equation}
for $s\in{\mathbb Z},\,i\in\{1,2\},\,\nu\in{\mathbb N}_0,$
\begin{equation}\label{eq:6ce}
\beta^{\ast(r)}_{\alpha,0,3+j}(z;\nu)=
\end{equation}
$$\sum\limits_{i=1}^2\left(\sum\limits_{k=0}^{\nu+\alpha}
(1-j+ij)\beta_{\alpha,0,i,k,\nu}^{(r)}
\left(\sum\limits_{\tau=1}^k
z^{k-\tau}(\tau)^{-i-j}\right)\right)=
$$
$$\sum\limits_{\sigma=0}^{\nu+\alpha-1}z^\sigma
\sum\limits_{\tau=1}^{\nu+\alpha-\sigma}
\sum\limits_{i=1}^2(1-j+ij)\beta_{\alpha,0,i,\sigma+\tau,\nu}^{(r)}
(\tau)^{-i-j}.$$
In view of (\ref{eq:6bh.0}) and (\ref{eq:6cd}), if
 $r=0,\,1,\,2,\,\alpha\in{\mathbb N}.$ then
\begin{equation}\label{eq:6cd1}
\beta^{\ast(r)}_{\alpha,0,1}(z;\nu)=
(z-1)\beta^{\ast\vee(r)}_{\alpha,0,1}(z;\nu),
\end{equation}
where $\beta^{\ast\vee(r)}_{\alpha,0,1}(z;\nu)\in{\mathbb Q}[z],$
when $\nu\in{\mathbb N}_0.$
In view of (\ref{eq:6bh.1}) and (\ref{eq:6cd}), 
\begin{equation}\label{eq:6cd2}
\beta^{\ast(3)}_{1,0,1}(z;\nu)=-(\nu+1)^2+
(z-1)\beta^{\ast\vee(3)}_{1,0,1}(z;\nu),
\end{equation}
where $\beta^{\ast\vee(3)}_{1,0,1}(z;\nu)\in{\mathbb Q}[z],$
when $\nu\in{\mathbb N}_0.$
In view of (\ref{eq:6bg.1.2}) -- (\ref{eq:6bg.3.1}), (\ref{eq:6cd}), 
\begin{equation}\label{eq:6bg.1.2a}
\beta_{\alpha,0,2}^{\ast(1)}(z;\nu)=\delta\beta_{\alpha,0,2}^{\ast(0)}(z;\nu),
\end{equation}
\begin{equation}\label{eq:6bg.1.1a}
\beta_{\alpha,0,1}^{\ast(1)}(z;\nu)=\delta\beta_{\alpha,0,1}^{\ast(0)}(z;\nu)-
\beta_{\alpha,0,2}^{\ast(0)}(z;\nu),
\end{equation}
\begin{equation}\label{eq:6bg.2.2a}
\beta_{\alpha,0,2}^{\ast(2)}(z;\nu)=
\delta^2\beta_{\alpha,0,2}^{\ast(0)}(z;\nu),
\end{equation}
\begin{equation}\label{eq:6bg.2.1a}
\beta_{\alpha,0,1}^{\ast(2)}(z;\nu)=
\delta^2\beta_{\alpha,0,1}^{\ast(0)}(z;\nu)-
2\delta\beta_{\alpha,0,2}^{\ast(0)}(z;\nu),
\end{equation}
\begin{equation}\label{eq:6bg.3.2a}
\beta_{\alpha,0,2}^{\ast(3)}(z;\nu)=
\delta^3\beta_{\alpha,0,2}^{\ast(0)}(z;\nu),
\end{equation}
\begin{equation}\label{eq:6bg.3.1a}
\beta_{\alpha,0,1}^{\ast(3)}(z;\nu)=
\delta^3\beta_{\alpha,0,1}^{\ast(0)}(z;\nu)-
3\delta^2\beta_{\alpha,0,2}^{\ast(0)}(z;\nu),
\end{equation}
where $\alpha\in{\mathbb N}.$ Clearly,
\begin{equation}\label{eq:6cg}
(-\delta)^kL_n(1/z)=L_{n-k}(1/z),
\end{equation}
where $k\in[0,+\infty)\cap{\mathbb Z},\,n\in{\mathbb Z},\,\vert z\vert>1,$
\begin{equation}\label{eq:6ch0}
L_1(1/z)=-\log(1-1/z),\,-\delta L_1(1/z)=
\end{equation}
$$L_0(1/z)=\frac1{z-1},\,\delta^2L_1(1/z)=$$
$$L_{-1}(1/z)=\frac1{z-1}+\frac1{(z-1)^2},\,-\delta^3L_1(1/z)=$$
$$L_{-2}(1/z)=\frac1{z-1}+\frac3{(z-1)^2}+\frac2{(z-1)^3}.$$ 
We apply the operator $\delta$ to the equality (\ref{eq:15.4.hr1f})
 for $r=0,\,1,\,2,\,\alpha\in{\mathbb N}.$ Then,
 in view of  (\ref{eq:6cg}), we obtain the eqvality
\begin{equation}\label{eq:15.4.hr2f}
\delta^{r+1}f_{\alpha,0,2+j}^{\ast}(z, \nu)-
j(\log(z))\delta^{r+1}f_{\alpha,0,2}^{\ast}(z, \nu)=
\end{equation}
$$j\delta^rf_{\alpha,0,2}^{\ast}(z, \nu)+$$
$$\left(\sum\limits_{i=1}^2
((1-j+ij)\delta\beta^{\ast(r)}_{\alpha,0,i}(z;\nu))L_{i+j}(1/z)\right)-
\delta\beta^{(r)}_{\alpha,0,3+j}(z;\nu)-$$
$$\left(\sum\limits_{i=1}^2
(1-j+ij)\beta^{\ast(r)}_{\alpha,0,i}(z;\nu)L_{i+j-1}(1/z)\right)=$$
$$j\left(\left(\sum\limits_{i=1}^2
\beta^{\ast(r)}_{\alpha,0,i}(z;\nu)L_i(1/z)\right)-
\beta^{(r)}_{\alpha,0,3}(z;\nu)\right)+$$
$$\left(\sum\limits_{i=1}^2
((1-j+ij)\delta\beta^{\ast(r)}_{\alpha,0,i}(z;\nu))L_{i+j}(1/z)\right)-
\delta\beta^{(r)}_{\alpha,0,3+j}(z;\nu)-$$
$$\left(\sum\limits_{i=1}^2
(1-j+ij)\beta^{\ast(r)}_{\alpha,0,i}(z;\nu)L_{i+j-1}(1/z)\right).$$
It follws from (\ref{eq:15.4.hr2f}) with $j=0$ that
\begin{equation}\label{eq:15.4.hr2f0}
\delta^{r+1}f_{\alpha,0,2}^{\ast}(z, \nu)=
\end{equation}
$$\left(\sum\limits_{i=1}^2
(\delta\beta^{\ast(r)}_{\alpha,0,i}(z;\nu))L_{i}(1/z)-
\beta^{\ast(r)}_{\alpha,0,i}(z;\nu)L_{i-1}(1/z)\right)-
\delta\beta^{(r)}_{\alpha,0,3}(z;\nu)=$$
$$(\delta\beta^{\ast(r)}_{\alpha,0,2}(z;\nu))L_2(1/z)+
(\delta\beta^{\ast(r)}_{\alpha,0,1}(z;\nu)-
\beta^{\ast(r)}_{\alpha,0,2}(z;\nu))L_{1}(1/z)-$$ 
$$\delta\beta^{\ast(r)}_{\alpha,0,3}(z;\nu)-
\beta^{\ast(r)}_{\alpha,0,1}(z;\nu)L_0(1/z).$$
In view of (\ref{eq:15.4.hr1f}) with $j=0,$ (\ref{eq:15.4.hr2f0}),
 (\ref{eq:6cd1}),
\begin{equation}\label{eq:15.4.hr3}
\beta^{\ast(r)}_{\alpha,0,2}(z;\nu)=
\delta\beta^{\ast(r-1)}_{\alpha,0,2}(z;\nu))=
\delta^r\beta^{\ast(0)}_{\alpha,0,2}(z;\nu),
\end{equation}
\begin{equation}\label{eq:15.4.hr4}
\beta^{\ast(r)}_{\alpha,0,1}(z;\nu))=
\delta\beta^{\ast(r-1)}_{\alpha,0,1}(z;\nu)-
\beta^{\ast(r-1)}_{\alpha,0,2}(z;\nu)=
\end{equation}
$$\delta^r\beta^{\ast(0)}_{\alpha,0,1}(z;\nu)-
r\delta^{r-1}\beta^{\ast(0)}_{\alpha,0,2}(z;\nu),$$
\begin{equation}\label{eq:15.4.hr5}
\beta^{\ast(r)}_{\alpha,0,3}(z;\nu)=
\delta\beta^{\ast(r-1)}_{\alpha,0,3}(z;\nu)+
\beta^{\ast(r-1)}_{\alpha,0,1}(z;\nu))L_0(1/z)=
\end{equation}
$$\delta\beta^{\ast(r-1)}_{\alpha,0,3}(z;\nu)+
\beta^{\ast\vee(r-1)}_{\alpha,0,1}(z;\nu),
$$
where $r=1,\,2,\,3,$ and $\alpha\in{\mathbb N}.$
The equalities (\ref{eq:6bg.1.2a}) -- (\ref{eq:6bg.3.1a}) 
follow from the equalities (\ref{eq:15.4.hr3}) and (\ref{eq:15.4.hr4}) again.
In view of (\ref{eq:6bi2}), (\ref{eq:6cd}), (\ref{eq:b})
 and (\ref{eq:15.4.hr3})  
\begin{equation}\label{eq:15.4.hr8}
\beta^{\ast(r)}_{\alpha,0,2}(z;\nu)=
\delta^rf^{\ast}_{\alpha,0,1}(z;\nu)\in{\mathbb N}[z],
\end{equation}
where $\alpha\in{\mathbb N},\,\nu\in{\mathbb N}_0,\,r=0,\,1,\,2,\,3.$
It follws from (\ref{eq:15.4.hr2f}) with $j=1$ that
\begin{equation}\label{eq:15.4.hr2f1}
\delta^{r+1}f_{\alpha,0,3}^{\ast}(z, \nu)=
(\log(z))\delta^{r+1}f_{\alpha,0,2}^{\ast}(z, \nu)+
\end{equation}
$$\left(\left(\sum\limits_{i=1}^2
\beta^{\ast(r)}_{\alpha,0,i}(z;\nu)L_i(1/z)\right)-
\beta^{(r)}_{\alpha,0,3}(z;\nu)\right)+$$
$$\left(\sum\limits_{i=1}^2
i(\delta\beta^{\ast(r)}_{\alpha,0,i}(z;\nu))L_{i+1}(1/z)\right)-
\delta\beta^{(r)}_{\alpha,0,4}(z;\nu)-$$
$$\left(\sum\limits_{i=1}^2
i\beta^{\ast(r)}_{\alpha,0,i}(z;\nu)L_{i}(1/z)\right)=
(\log(z))\delta^{r+1}f_{\alpha,0,2}^{\ast}(z, \nu)+$$
$$\left(\sum\limits_{i=1}^2
i(\delta\beta^{\ast(r)}_{\alpha,0,i}(z;\nu))L_{i+1}(1/z)\right)-
\delta\beta^{(r)}_{\alpha,0,4}(z;\nu)-$$
$$\beta^{\ast(r)}_{\alpha,0,2}(z;\nu)L_{2}(1/z)-
(\beta^{(r)}_{\alpha,0,3}(z;\nu)+
\delta\beta^{(r)}_{\alpha,0,4}(z;\nu))=$$
$$(\log(z))\delta^{r+1}f_{\alpha,0,2}^{\ast}(z, \nu)+
2(\delta\beta^{\ast(r)}_{\alpha,0,2}(z;\nu))L_3(1/z)+$$
$$(\delta\beta^{\ast(r)}_{\alpha,0,1}(z;\nu)-
\beta^{\ast(r)}_{\alpha,0,2}(z;\nu))L_2(1/z)-
(\delta\beta^{(r)}_{\alpha,0,4}(z;\nu)+\beta^{\ast(r)}_{\alpha,0,3}(z;\nu)).$$
In view of (\ref{eq:15.4.hr1f}) with $j=1,$ (\ref{eq:15.4.hr2f1}),
\begin{equation}\label{eq:15.4.hr2f2}
\beta^{\ast(r+1)}_{\alpha,0,2}(z;\nu)=\delta\beta^{\ast(r)}_{\alpha,0,2}(z;\nu)=
\delta^{r+1}\beta^{\ast(0)}_{\alpha,0,2}(z;\nu),
\end{equation}
\begin{equation}\label{eq:15.4.hr2f3}
\beta^{\ast(r+1)}_{\alpha,0,1}(z;\nu)=
\delta\beta^{\ast(r)}_{\alpha,0,2}(z;\nu)-
\beta^{\ast(r)}_{\alpha,0,2}(z;\nu)=
\end{equation}
$$\delta^{r+1}\beta^{(0)}_{\alpha,0,1}(z;\nu)-
\delta^r\beta^{(0)}_{\alpha,0,2}(z;\nu),$$
where $r=0,\,1,\,2,$
and we obtain (\ref{eq:15.4.hr3}) -- (\ref{eq:15.4.hr3}) again. Moreover,
\begin{equation}\label{eq:15.4.hr2f4}
\beta^{\ast(r+1)}_{\alpha,0,4}(z;\nu)=
\delta\beta^{\ast(r)}_{\alpha,0,4}(z;\nu)+
\beta^{\ast(r)}_{\alpha,0,3}(z;\nu),
\end{equation}
where $r=0,\,1,\,2.$
If we take now $z\in(1,+\infty)$ and will tend $z$ to $1,$ then,
in view of (\ref{eq:15.4.hr1f}), (\ref{eq:6cd1}), (\ref{eq:6cd2})
 and (\ref{eq:6ch0})
\begin{equation}\label{eq:15.4.hr6}
\delta^rf_{\alpha,0,2+j}^{\ast}(1,\nu)
=\lim\limits_{z\to1+0}\delta^rf_{\alpha,0,2}^{\ast}(z, \nu)=
\end{equation}
$$(1-i+ij\beta^{\ast(r)}_{\alpha,0,2}(1;\nu)\zeta(2+j)-
\beta^{\ast(r)}_{\alpha,0,3+j}(1;\nu),$$
where $r=0,\,1,\,2,\,j=0,\,1,$
\begin{equation}\label{eq:15.4.hr7}
\lim\limits_{z\to1+0} (z-1)\delta^3f_{\alpha,0,2}^{\ast}(z, \nu)=0,
\end{equation}
if $\alpha\in{\mathbb N}.$
In view of (\ref{eq:6cd}), (\ref{eq:6ce}),
 (\ref{eq:r0a0i2k01n01}) -- (\ref{eq:r0a1i1k2n1}),
\begin{equation}\label{eq:a0r0i1n0}
\beta^{\ast(0)}_{0,0,1}(z;0)=\beta_{0,0,1,0,0}^{(0)}=0,
\end{equation}  
\begin{equation}\label{eq:a0r0i2n0}
\beta^{\ast(0)}_{0,0,2}(z;0)=\beta_{0,0,2,0,0}^{(0)}=1,
\end{equation}
\begin{equation}\label{eq:a0r0i34n0}
\beta^{\ast(0)}_{0,0,3}(z;0)=\beta^{\ast(0)}_{0,0,4}(z;0)=0,
\end{equation}
\begin{equation}\label{eq:a0r0i1n1}
\beta^{\ast(0)}_{0,0,1}(z;1)=\beta_{0,0,1,0,1}^{(0)}+
\beta_{0,0,1,1,1}^{(0)}z=-4+4z,
\end{equation}
\begin{equation}\label{eq:a0r0i2n1}
\beta^{\ast(0)}_{0,0,2}(z;1)=\beta_{0,0,2,0,1}^{(0)}+
\beta_{0,0,2,1,1}^{(0)}z=1+4z,
\end{equation}
\begin{equation}\label{eq:a0r0i3n1}
\beta^{\ast(0)}_{0,0,3}(z;1)=\beta_{0,0,1,1,1}^{(0)}+
\beta_{0,0,2,1,1}^{(0)}=8,
\end{equation}
\begin{equation}\label{eq:a0r0i4n1}
\beta^{\ast(0)}_{0,0,4}(z;1)=\beta_{0,0,1,1,1}^{(0)}+
2\beta_{0,0,2,1,1}^{(0)}=12,
\end{equation}
\begin{equation}\label{eq:a1r0i1n0}
\beta^{\ast(0)}_{1,0,1}(z;0)=\beta_{1,0,1,0,0}^{(0)}+
\beta_{1,0,1,1,0}^{(0)}z=-2+2z,
\end{equation}
\begin{equation}\label{eq:a1r0i2n0}
\beta^{\ast(0)}_{1,0,2}(z;0)=\beta_{1,0,2,0,0}^{(0)}+
\beta_{1,0,2,1,0}^{(0)}z=1+z,
\end{equation}
\begin{equation}\label{eq:a1r0i3n0}
\beta^{\ast(0)}_{1,0,3}(z;0)=\beta_{1,0,1,1,0}^{(0)}+
\beta_{1,0,2,1,0}^{(0)}=3,
\end{equation}
\begin{equation}\label{eq:a1r0i4n0}
\beta^{\ast(0)}_{1,0,4}(z;0)=\beta_{1,0,1,1,0}^{(0)}+
2\beta_{1,0,2,1,0}^{(0)}=4,
\end{equation}
\begin{equation}\label{eq:a1r0i1n1}
\beta^{\ast(0)}_{1,0,1}(z;1)=\beta_{1,0,1,0,1}^{(0)}+
\beta_{1,0,1,1,1}^{(0)}z+
\end{equation}
$$\beta_{1,0,1,2,1}^{(0)}z^2=-5-16z+21z^2=(z-1)(21z+5),$$
\begin{equation}\label{eq:a1r0i2n1}
\beta^{\ast(0)}_{1,0,2}(z;1)=\beta_{1,0,2,0,1}^{(0)}+
\beta_{1,0,2,1,1}^{(0)}z+
\end{equation}
$$\beta_{1,0,2,2,1}^{(0)}z^2=1+16z+9z^2,$$
\begin{equation}\label{eq:a1r0i3n1}
\beta^{\ast(0)}_{1,0,3}(z;1)=\beta_{1,0,1,1,1}^{(0)}+\beta_{1,0,2,1,1}^{(0)}+
\end{equation}
$$\frac12\beta_{1,0,1,2,1}^{(0)}+\frac14\beta_{1,0,2,2,1}^{(0)}+
(\beta_{1,0,1,2,1}^{(0)}+\beta_{1,0,2,2,1}^{(0)})z=$$
$$-16+16+\frac12\times21+\frac14\times9+(21+9)z=\frac{51}4+30z,$$
\begin{equation}\label{eq:a1r0i4n1}
\beta^{\ast(0)}_{1,0,4}(z;1)=\beta_{1,0,1,1,1}^{(0)}+2\beta_{1,0,2,1,1}^{(0)}+
\end{equation}
$$\frac14\beta_{1,0,1,2,1}^{(0)}+2\times\frac18\beta_{1,0,2,2,1}^{(0)}+
(\beta_{1,0,1,2,1}^{(0)}+2\beta_{1,0,2,2,1}^{(0)})z=$$
$$-16+2\times16+\frac14\times21+\frac14\times9+(21+18)z=\frac{47}2+39z.$$
\begin{equation}\label{eq:a1r0i1n2}
\beta^{\ast(0)}_{1,0,1}(z;2)=\beta_{1,0,1,0,2}^{(0)}+\beta_{1,0,1,1,2}^{(0)}z+
\end{equation}
$$\beta_{1,0,1,2,2}^{(0)}z^2+\beta_{1,0,1,3,2}^{(0)}z^3=$$
$$-\frac{20}3-216z-54z^2+\frac{830}3z^3=(z-1)(830z^2+668z+20)/3.$$
\begin{equation}\label{eq:a1r0i2n2}
\beta^{\ast(0)}_{1,0,2}(z;2)=\beta_{1,0,2,0,2}^{(0)}+\beta_{1,0,2,1,2}^{(0)}z+
\beta_{1,0,2,2,2}^{(0)}z^2+
\end{equation}
$$\beta_{1,0,2,3,2}^{(0)}z^3=1+81z+324z^2+100z^3.$$
\begin{equation}\label{eq:a1r0i3n2}
\beta^{\ast(0)}_{1,0,3}(z;2)=\beta_{1,0,1,1,2}^{(0)}+\beta_{1,0,2,1,2}^{(0)}+
\end{equation}
$$\frac12\beta_{1,0,1,2,2}^{(0)}+\frac14\beta_{1,0,2,2,2}^{(0)}+
\frac13\beta_{1,0,1,3,2}^{(0)}+\frac19\beta_{1,0,2,3,2}^{(0)}+$$
$$(\beta_{1,0,1,2,2}^{(0)}+\beta_{1,0,2,2,2}^{(0)})z+
\left(\frac12\beta_{1,0,1,3,2}^{(0)}+\frac14\beta_{1,0,2,3,2}^{(0)}\right)z+
$$
$$(\beta_{1,0,1,3,2}^{(0)}+\beta_{1,0,2,3,2}^{(0)})z^2=$$
$$-216+81-54/2+324/4+(830+100)/9+((-54+324-54+415/3+100/4)z+$$
$$(830/3+100)z^2=\frac{67}3+\frac{1300}3z+\frac{1130}3z^2.$$
\begin{equation}\label{eq:a1r0i4n2}
\beta^{\ast(0)}_{1,0,4}(z;2)=\beta_{1,0,1,1,2}^{(0)}+2\beta_{1,0,2,1,2}^{(0)}+
\end{equation}
$$\frac14\beta_{1,0,1,2,2}^{(0)}+2\times\frac18\beta_{1,0,2,2,2}^{(0)}+$$
$$\frac19\beta_{1,0,1,3,2}^{(0)}+2\times\frac1{27}\beta_{1,0,2,3,2}^{(0)}+$$
$$(\beta_{1,0,1,2,2}^{(0)}+2\beta_{1,0,2,2,2}^{(0)})z+
\left(\frac14\beta_{1,0,1,3,2}^{(0)}+2\times\frac18\beta_{1,0,2,3,2}^{(0)}
\right)z+
$$
$$
(\beta_{1,0,1,3,2}^{(0)}+2\beta_{1,0,2,3,2}^{(0)})z^2=
$$
$$-216+162+(-54+324)/4+(830+200)/27+((-57+648+(830/3+100)/4)z+$$
$$(830/3+200)z^2=\frac{2789}{54}+\frac{4129}6z+\frac{1430}3z^2.$$
In view of (\ref{eq:15.4.hr3}) -- (\ref{eq:15.4.hr5}), (\ref{eq:15.4.hr2f4})
and (\ref{eq:a1r0i1n0}) -- (\ref{eq:a1r0i4n2}),
\begin{equation}\label{eq:a1r1i1n0}
\beta^{\ast(1)}_{1,0,1}(z;0)=\delta\beta^{\ast(0)}_{1,0,1}(z;0)-
\beta^{\ast(0)}_{1,0,2}(z;0)=
\end{equation}
$$2z-1-z=z-1,$$
\begin{equation}\label{eq:a1r1i2n0}
\beta^{\ast(1)}_{1,0,2}(z;0)=\delta\beta^{\ast(0)}_{1,0,2}(z;0)=z,
\end{equation}
\begin{equation}\label{eq:a1r1i3n0}
\beta^{\ast(1)}_{1,0,3}(z;0)=\delta\beta^{\ast(0)}_{1,0,3}(z;0)+
\beta^{\ast\vee(0)}_{1,0,1}(z;0)=2
\end{equation}
\begin{equation}\label{eq:a1r1i4n0}
\beta^{\ast(1)}_{1,0,4}(z;0)=\delta\beta^{\ast(0)}_{1,0,4}(z;0)+
\beta^{\ast(0)}_{1,0,3}(z;0)=3,
\end{equation}
\begin{equation}\label{eq:a1r1i1n1}
\beta^{\ast(1)}_{1,0,1}(z;1)=\delta\beta^{\ast(0)}_{1,0,1}(z;1)-
\beta^{\ast(0)}_{1,0,2}(z;1)=-16z+
\end{equation}
$$42z^2-(1+16z+9z^2)=-1-32z+33z^2=(z-1)(33z+1),$$
\begin{equation}\label{eq:a1r1i2n1}
\beta^{\ast(1)}_{1,0,2}(z;1)=\delta\beta^{\ast(0)}_{1,0,2}(z;1)=
16z+18z^2,
\end{equation}
\begin{equation}\label{eq:a1r1i3n1}
\beta^{\ast(1)}_{1,0,3}(z;1)=\delta\beta^{\ast(0)}_{1,0,3}(z;1)+
\beta^{\ast\vee(0)}_{1,0,1}(z;1)=
\end{equation}
$$30z+21z+5=51z+5,$$
\begin{equation}\label{eq:a1r1i4n1}
\beta^{\ast(1)}_{1,0,4}(z;1)=\delta\beta^{\ast(0)}_{1,0,4}(z;1)+
\beta^{\ast(0)}_{1,0,3}(z;1)=
\end{equation}
$$39z+\frac{51}4+30z=69z+\frac{51}4,$$
\begin{equation}\label{eq:a1r1i1n2}
\beta^{\ast(1)}_{1,0,1}(z;2)=\delta\beta^{\ast(0)}_{1,0,1}(z;2)-
\beta^{\ast(0)}_{1,0,2}(z;2)=
\end{equation}
$$-216z-108z^2+830z^3-(1+81z+324z^2+100z^3)=$$
$$-1-297z-432z^2+730z^3=(z-1)(730z^2+298z+1),$$
\begin{equation}\label{eq:a1r1i2n2}
\beta^{\ast(1)}_{1,0,2}(z;2)=\delta\beta^{\ast(0)}_{1,0,2}(z;2)=
81z+648z^2+300z^3,
\end{equation}
\begin{equation}\label{eq:a1r1i3n2}
\beta^{\ast(1)}_{1,0,3}(z;2)=\delta\beta^{\ast(0)}_{1,0,3}(z;2)+
\beta^{\ast\vee(0)}_{1,0,1}(z;2)=
\end{equation}
$$\frac{1300}3z+\frac{2260}3z^2+\frac{830}3z^2+\frac{668}3z+\frac{20}3=$$
$$1030z^2+656z+\frac{20}3,$$
\begin{equation}\label{eq:a1r1i4n2}
\beta^{\ast(1)}_{1,0,4}(z;2)=\delta\beta^{\ast(0)}_{1,0,4}(z;2)+
\beta^{\ast(0)}_{1,0,3}(z;2)=
\end{equation}
$$\frac{4129}6z+\frac{2860}3z^2+\frac{67}3+\frac{1300}3z+\frac{1130}3z^2=$$
$$\frac{67}3+\frac{6729}6z+1330z^2,$$
In view of (\ref{eq:15.4.hr3}) -- (\ref{eq:15.4.hr5}), (\ref{eq:15.4.hr2f4})
and (\ref{eq:a1r1i1n0}) -- (\ref{eq:a1r1i4n2}),
\begin{equation}\label{eq:a1r2i1n0}
\beta^{\ast(2)}_{1,0,1}(z;0)=\delta\beta^{\ast(1)}_{1,0,1}(z;0)-
\beta^{\ast(1)}_{1,0,2}(z;0)=z-z=0,
\end{equation}
\begin{equation}\label{eq:a1r2i2n0}
\beta^{\ast(2)}_{1,0,2}(z;0)=\delta\beta^{\ast(1)}_{1,0,2}(z;0)=z,
\end{equation}
\begin{equation}\label{eq:a1r2i3n0}
\beta^{\ast(2)}_{1,0,3}(z;0)=\delta\beta^{\ast(1)}_{1,0,3}(z;0)+
\beta^{\ast\vee(1)}_{1,0,1}(z;0)=1,
\end{equation}
\begin{equation}\label{eq:a1r2i4n0}
\beta^{\ast(2)}_{1,0,4}(z;0)=\delta\beta^{\ast(1)}_{1,0,4}(z;0)+
\beta^{\ast(1)}_{1,0,3}(z;0)=2
\end{equation}
\begin{equation}\label{eq:a1r2i1n1}
\beta^{\ast(2)}_{1,0,1}(z;1)=\delta\beta^{\ast(1)}_{1,0,1}(z;1)-
\beta^{\ast(1)}_{1,0,2}(z;1)=
\end{equation}
$$-32z+66z^2-(16z+18z^2)=-48z+48z^2=48z(z-1).$$
\begin{equation}\label{eq:a1r2i2n1}
\beta^{\ast(2)}_{1,0,2}(z;1)=\delta\beta^{\ast(1)}_{1,0,2}(z;1)=
16z+36z^2,
\end{equation}
\begin{equation}\label{eq:a1r2i3n1}
\beta^{\ast(2)}_{1,0,3}(z;1)=\delta\beta^{\ast(1)}_{1,0,3}(z;1)+
\beta^{\ast\vee(1)}_{1,0,1}(z;1)=
\end{equation}
$$51z+33z+1=84z+1,$$
\begin{equation}\label{eq:a1r2i4n1}
\beta^{\ast(2)}_{1,0,4}(z;1)=\delta\beta^{\ast(1)}_{1,0,4}(z;1)+
\beta^{\ast(1)}_{1,0,3}(z;1)=
\end{equation}
$$69z+51z+5=120z+5,$$
\begin{equation}\label{eq:a1r2i1n2}
\beta^{\ast(2)}_{1,0,1}(z;2)=\delta\beta^{\ast(1)}_{1,0,1}(z;2)-
\beta^{\ast(1)}_{1,0,2}(z;2)=
\end{equation}
$$-297z-864z^2+2190z^3-(81z+648z^2+300z^3)=$$
$$z(-378-1512z+1890z^2)=378(z-1)z(5z+1),$$
\begin{equation}\label{eq:a1r2i2n2}
\beta^{\ast(2)}_{1,0,2}(z;2)=\delta\beta^{\ast(1)}_{1,0,2}(z;2)=
81z+1296z^2+900z^3,
\end{equation}
\begin{equation}\label{eq:a1r2i3n2}
\beta^{\ast(0)}_{1,0,3}(z;2)=\delta\beta^{\ast(1)}_{1,0,3}(z;2)+
\beta^{\ast\vee(1)}_{1,0,1}(z;2)
\end{equation}
$$2060z^2+656z+730z^2+298z+1=2790z^2+954z+1,$$
\begin{equation}\label{eq:a1r2i4n2}
\beta^{\ast(2)}_{1,0,4}(z;2)=\delta\beta^{\ast(1)}_{1,0,4}(z;2)+
\beta^{\ast(1)}_{1,0,3}(z;2)=
\end{equation}
$$\frac{6729}6z+2660z^2+1030z^2+656z+\frac{20}3=3690z^2+\frac{3555}2z+
\frac{20}3,$$
{\begin{center}\large\bf
\S 7. End of the proof of theorem A.
 \end{center}}
Let $y^{(\kappa)}_{1,0,i,k}(z;\nu)$ denotes $i-$th element of
 the column $Y^{(\kappa)}_{1,0,,k}(z;\nu)$
in (\ref{eq:15.1.h}). Then, in view of (\ref{eq:15e}), (\ref{eq:15.1.b}),
 (\ref{eq:15.1.h}),
\begin{equation}\label{eq:7.1}
y^{(\kappa)}_{1,0,i+1,k}(z,\nu)=\delta^i f^{(\kappa)}_{\alpha,0,k}(z,\nu)
\end{equation}
for
$\kappa=0,\,1,\,i=1,\,2,\,3,\,k=1,\,2,\,3,\,\vert z\vert>1,
\,\nu\in{\mathbb N}_0.$
We denote $v^{\ast\kappa\ast}_{1,0,i,j}(\nu)$ the expression, which 
stands in the matrix $V^{\ast\kappa\ast}_{1,0}(\nu)$ in intersection
of $i$-th row and $j$-th column, where $i=1,\,2,\,3,\,4,\,j=1,\,2,\,3,\,4.$ 
Let
\begin{equation}\label{eq:93bd1}
 D_{\alpha,0}(z,\nu,w)=z(w^2+w(1-\alpha)-\mu_\alpha)^{2}-w^{4},
\end{equation}
In view of (\ref{eq:93bd})
\begin{equation}\label{eq:93bd1a}
\frac1z D_{\alpha,0}(z,\nu,w)=(1-1/z)w^4+
\sum\limits_{k=0}^3r_{\alpha,0,k+1}(\nu)w^k,
\end{equation}
It follows from general properties of Mejer's functions that
\begin{equation}\label{eq:93bd2}
 D_{\alpha,0}(z,\nu,\delta)f_{\alpha,0,k}(z,\nu)=0,
\end{equation}
where
$\vert z\vert>1,-3\pi/2<\arg(z)\le\pi/2,\log(z)=\ln(\vert z\vert)+i\arg(z),\,
k=1,\,2,\,3.$ Therefore,  in view of (\ref{eq:15e}), (\ref{eq:15.1.b}),
 (\ref{eq:15.1.h}),
\begin{equation}\label{eq:93bd3}
y^{(\kappa)}_{2\kappa+1,0,k}(z,\nu)=
-(1-1/z)\delta^4 f_{\alpha,0,k}(z,\nu)
\end{equation}
where
$$\kappa=0,\,1,,\,\vert z\vert>1,-3\pi/2<\arg(z)\le\pi/2,$$
$$\log(z)=\ln(\vert z\vert)+i\arg(z),\,k=1,\,2,\,3.$$
 In view of (\ref{eq:b}) -- (\ref{eq:ac}), (\ref{eq:acrf}), 
\begin{equation}\label{eq:7.2}
\lim\limits_{z\to1+0}(z-1)\delta^4 f_{1,0,2}(z,\nu)=
\end{equation}
$$\lim\limits_{z\to1+0}(z-1)
\left(O(1)\ln\left(1-\frac1z\right)+\frac1{z-1}\right)=1,$$
\begin{equation}\label{eq:7.3}
\lim\limits_{z\to1+0}(z-1)\delta^4 f_{1,0,k}(z,\nu)=0,
\end{equation}
if $k-2=\pm1,$
\begin{equation}\label{eq:7.4}
\lim\limits_{z\to1+0}(\log(z))\delta^i f_{1,0,k}(z,\nu)=
\lim\limits_{z\to1+0}(z-1)\delta^i f_{1,0,k}(z,\nu)=0,
\end{equation}
if $i=0,\,1,\,2,\,3,k=1,\,2,\,3.$
Hence, if we tend $z\in(1,+\infty)$ to $1,$ then,
 in view of (\ref{eq:15e}), (\ref{eq:15.1.b}), we obtain
the equalities
\begin{equation}\label{eq:7.5}
y^{(\kappa)}_{1,0,2\kappa+1,1}(1,\nu)=y_{1,0,2\kappa+1,3}(1,\nu)=0,
y_{1,0,2\kappa+1,2}(1,\nu)=-1.
\end{equation}
In view of (\ref{eq:15.1.aj}), (\ref{eq:15.2.a11}) -- (\ref{eq:1as43}),
(\ref{eq:7.1}), (\ref{eq:93bd3}),
\begin{equation}\label{eq:7.6}
-a_{1,0,i,2\kappa+1}^{\ast\kappa\ast}(1;\nu)(1-1/z)\delta^4 f_{1,0,k}(z,\nu)+
\end{equation}
$$\left(\sum\limits_{j=1}^2
a_{1,0,i,j+1-\kappa}^{\ast\kappa\ast}(1;\nu)
\delta^{j-\kappa}f_{1,0,k}(z,\nu)\right)-
$$
$$
(z-1)v_{1,0,i,1}^{\ast\kappa\ast}(\nu)(1-1/z)\delta^4 f_{1,0,k}(z,\nu)+
$$
$$(z-1)\sum\limits_{j=1}^2
v_{1,0,i+1,j+1-\kappa}^{\ast\ast}(\nu)\delta^{j-\kappa}f_{1,0,k}(z,\nu)=
$$
$$(-2)^\kappa\mu_1(\nu)^{2-\kappa}\nu^5\delta^{i-1}f_{\alpha,0,k}(z,\nu-1),$$
where
$i=3-2\kappa,\,2,\,k=1,\,2,\,3,\,\vert z\vert>1,\,-3\pi/2<\arg(z)\le\pi/2$ and
$\nu$ run over the set $M_1^{\ast}=((-\infty,-2]\cup[1,+\infty))\cap{\mathbb Z}.$
We tend $z\in(1,+\infty)$ to $1$ now and obtain the equalities
\begin{equation}\label{eq:7.7}
a_{1,0,i,1}^{\ast\kappa\ast}(1;\nu)(k-1)(k-3)+
\end{equation}
$$\left(\sum\limits_{j=1}^2
a_{1,0,i,j+1-\kappa}^{\ast\kappa\ast}(1;\nu)
(\delta^{j-\kappa}f_{1,0,k})(1,\nu)\right)=$$
$$(-2)^\kappa\mu_1(\nu)^{2-\kappa}\nu^5\delta^{i-1}if_{1,0,k}(1,\nu-1),$$
where
$i=2-\kappa,\,3-\kappa,\,k=1,\,2,\,3$  and
$\nu\in M_1^{\ast}=((-\infty,-2]\cup[1,+\infty))\cap{\mathbb Z}.$
Replacing in (\ref{eq:7.7}) $\nu\in M_1^\ast$ by 
$$\nu:=-\nu-2\in M_1^{\ast\ast}=((-\infty,-3]\cup[0,+\infty))\cap{\mathbb Z},$$
 and taking in account the equality (\ref{eq:15.1.i}) we obtain the equalities
\begin{equation}\label{eq:7.8}
a_{1,0,i,1}^{\ast\kappa\ast}(1;-\nu-2)(k-1)(k-3)+
\end{equation}
$$\left(\sum\limits_{j=1}^2
a_{1,0,i,j+1-\kappa}^{\ast\kappa\ast}(1;-\nu-2)
(\delta^{j-\kappa}f_{1,0,k})(1,\nu)\right)=$$
$$(-2)^\kappa\mu_1(\nu)^{2-\kappa}\nu^5\delta^{i-1}f_{1,0,k}(1,\nu+1),$$
where
$i=2-\kappa,\,3-\kappa,\,k=1,\,2,\,3$  and
$\nu\in M_1^{\ast\ast}=((-\infty,-3]\cup[0,+\infty))\cap{\mathbb Z}.$ Let
\begin{equation}\label{eq:7.12}
\vec w^{(\kappa)}_{i,j}(\nu)=
\left(\matrix a_{1,0,i+1-\kappa,j+1-\kappa}^{\ast\ast}(1;-\nu-2)\\
                                   (1+(-1)^{i+j})/2\\
a_{1,0,i+1-\kappa,j+1-\kappa}^{\ast\ast}(1;\nu)
\endmatrix\right),\end{equation}
\begin{equation}\label{eq:7.13}
W^{(\kappa)}_i(\nu)=\left(\matrix \vec w_{i,1}(\nu)&\vec w_{i,2}(\nu)
\endmatrix\right),
\end{equation}
$$Y^{\ast\kappa\ast\ast}_k(\nu)=
\left(\matrix (\delta^{1-\kappa} f_{1,0,k})(1,\nu)\\
(\delta^{2-\kappa}f_{1,0,k})(1,\nu)\endmatrix\right),$$
\begin{equation}\label{eq:7.14}
Y^{\ast\kappa\ast\ast\ast}_{i,k}(\nu)=
\left(\matrix 
-(-2)^\kappa\mu_1(\nu)^{2-\kappa}(\nu+2)^5\delta^{i-\kappa}f_{1,0,k}(1,\nu+1)\\
(\delta^{i-\kappa} f_{1,0,k})(1,\nu)\\
(-2)^\kappa\mu_1(\nu)^{2-\kappa}\nu^5\delta^{i-\kappa}f_{\alpha,0,k}(1,\nu-1)
\endmatrix\right),\end{equation}
where $\kappa=0,\,1,\,k=1,\,3,\,i=1,\,2,\, 
\nu\in M_1^{\ast\ast\ast\ast}=((-\infty,-3]\cup[1,+\infty))\cap{\mathbb Z}.$
Let further
\begin{equation}\label{eq:7.11}
\vec w^{(\kappa)}_{i,3}(\nu)=
\left(\matrix w^{(\kappa)}_{i,3,1}(\nu)\\
                w^{(\kappa)}_{i,3,2}(\nu)\\
                 w^{(\kappa)}_{i,3,3}(\nu)
\endmatrix\right)=[\vec w^{(\kappa)}_{i,1}(\nu),\vec w^{(\kappa)}_{i,2}(\nu)]
\end{equation}
is vector product of $\vec w^{(\kappa)}_{i,1}(\nu)$
 and $\vec w^{(\kappa)}_{i,2}(\nu),$
and let $\bar w^{(\kappa)}_{3,i}(\nu)=(\vec w^{(\kappa)}_{3,i}(\nu))^t$ is row
conyugte to the column $\vec w(\kappa)_{3,i}(\nu).$
Then for scalar products
 $(\vec w^{(\kappa)}_{i,3}(\nu),\vec w^{(\kappa)}_{i,j}(\nu))$
we have the equalities
$$\bar w^{(\kappa)}_{i,3}(\nu)\vec w^{\kappa}_{i,j}(\nu)=
(\vec w^{(\kappa)}_{i,3}(\nu),\vec w^{(\kappa)}_{i,j}(\nu))=0,$$
where $\kappa=0,\,1,\,i=1,\,2,\,j=1,\,2$
and
$$
\nu\in M_1^{\ast\ast\ast\ast}=((-\infty,-3]\cup[1,+\infty))\cap{\mathbb Z}.
$$
Therefore
\begin{equation}\label{eq:7.9}
\bar w^{(\kappa)}_{i,3}(\nu)W^{(\kappa)}_i(\nu)=
\left(\matrix 0&0\endmatrix\right),
\end{equation}
where
$\kappa=0,\,1,\,i=1,\,2$ and
$\nu\in M_1^{\ast\ast\ast\ast}=((-\infty,-3]\cup[1,+\infty))\cap{\mathbb Z}.$

In view of (\ref{eq:7.7}) (\ref{eq:7.8}) and (\ref{eq:7.9}),
\begin{equation}\label{eq:7.10}
\bar w(\kappa)_{i,3}(\nu)Y^{\ast\kappa\ast\ast\ast}_{i,k}(\nu)=
\bar w(\kappa)_{i,3}(\nu)W(\kappa)_i(\nu)Y^{\ast\kappa\ast\ast}_k(\nu)=0,
\end{equation}
where
$\kappa=0,\,1,\,i=1,\,2,\,k=1,\,3$ and
$\nu\in M_1^{\ast\ast\ast\ast}=((-\infty,-3]\cup[1,+\infty))\cap{\mathbb Z}.$
In view of (\ref{eq:7.11}) -- (\ref{eq:7.13}),
 (\ref{eq:15.2.e}), and (\ref{eq:15.2.a11}) -- (\ref{eq:1as43})
\begin{equation}\label{eq:7w0131} 
w^(0)_{1,3,1}(\nu)=a_{1,0,2,3}^{\ast0\ast}(1;\nu)=
-3\tau^4(\tau-1)(2\tau-1)^3, 
\end{equation}
\begin{equation}\label{eq:7w1131} 
w^(1)_{1,3,1}(\nu)=a_{1,0,1,2}^{\ast1\ast}(1;\nu)=
\end{equation}
$$8\tau^3(2\tau-1)(\tau^3-(\tau-1)^3),$$

\begin{equation}\label{eq:7w0231} 
w^{(0)}_{2,3,1,}(\nu)=-a_{1,0,3,2}^{\ast0\ast}(1;\nu)=
\end{equation}
$$2\tau^5(\tau-1)^2(2\tau-1)(\tau^3-(\tau-1)^3),$$
\begin{equation}\label{eq:7w1231} 
w^(1)_{2,3,1}(\nu)=-a_{1,0,2,1}^{\ast1\ast}(1;\nu)=
a_{1,0,2,3}^{\ast0\ast}(1;\nu)=
\end{equation}
$$-3\tau^4(\tau-1)(2\tau-1)^3=w^(0)_{1,3,1}(\nu).$$
Since $\tau_\alpha(-\nu-1-\alpha)=-\tau_\alpha(\nu),$
it follows from (\ref{eq:7w0131}) -- (\ref{eq:7w1231}) that  
\begin{equation}\label{eq:7w0133} 
w^{(0)}_{1,3,3}(\nu)=-a_{1,0,2,3}^{\ast0\ast}(1;-\nu-2)=
3\tau^4(\tau+1)(2\tau+1)^3,
\end{equation}
\begin{equation}\label{eq:7w1133} 
w^{(0)}_{1,3,3}(\nu)=-a_{1,0,1,2}^{\ast1\ast}(1;-\nu-2)=
\end{equation}
$$-8\tau^3(2\tau+1)((\tau+1)^3-\tau^3),$$
\begin{equation}\label{eq:7w0233} 
w^{(0)}_{2,3,3}(\nu)=a_{1,0,3,2}^{\ast0\ast}(1;-\nu-2)= 
\end{equation}
$$-2\tau^5(\tau+1)^2(2\tau+1)((\tau+1)^3-\tau^3),$$
\begin{equation}\label{eq:7w1233} 
w^(1)_{2,3,3}(\nu)=a_{1,0,2,1}^{\ast1\ast}(1;-\nu-2)=
-a_{1,0,2,3}^{\ast0\ast}(1;-\nu-2)=
\end{equation}
$$3\tau^4(\tau+1)(2\tau+1)^3=w^(0)_{1,3,3}(\nu).$$
Further we have
\begin{equation}\label{eq:7w0132} 
w^{(0)}_{1,3,2}(\nu)=
\end{equation}
$$-\det\left(\matrix
a_{1,0,2,2}^{\ast(0)\ast}(1;-\nu-2)&a_{1,0,2,3}^{\ast(0)\ast}(1;-\nu-2)\\
a_{1,0,2,2}^{\ast0\ast}(1;\nu)&a_{1,0,2,3}^{\ast0\ast}(1;\nu)
\endmatrix\right)=$$ 
$$a_{1,0,2,2}^{\ast\ast}(1;\nu)a_{1,0,2,3}^{\ast\ast}(1;-\nu-2)-
a_{1,0,2,3}^{\ast\ast}(1;\nu)a_{1,0,2,2}^{\ast\ast}(1;-\nu-2)=$$
$$\tau^5(\tau-1)(\tau^3+2(2\tau-1)^3)
(-3\tau^4(\tau+1)(2\tau+1)^3)-$$
$$(-3\tau^4(\tau-1)(2\tau-1)^3)
(-\tau^5(\tau+1)(\tau^3+2(2\tau+1)^3))=$$
$$-3\tau^9(\tau^2-1)
(\tau^3((2\tau-1)^3+(2\tau+1)^3)+4(4\tau^2-1)^3)=$$
$$-12\tau^9(\tau^2-1)(68\tau^6-45\tau^4+12\tau^2-1),$$
\begin{equation}\label{eq:7w0232} 
w^{(0)}_{3,2,2}(\nu)=
\end{equation}
$$-\det\left(\matrix
a_{1,0,3,2}^{\ast\ast}(1;-\nu-2)&a_{1,0,3,3}^{\ast\ast}(1;-\nu-2)\\
a_{1,0,3,2}^{\ast\ast}(1;\nu)&a_{1,0,3,3}^{\ast\ast}(1;\nu)
\endmatrix\right)=$$ 
$$a_{1,0,3,2}^{\ast\ast}(1;\nu)a_{1,0,3,3}^{\ast\ast}(1;-\nu-2)-
a_{1,0,3,3}^{\ast\ast}(1;\nu)a_{1,0,3,2}^{\ast\ast}(1;-\nu-2)=$$
$$-2\tau^5(\tau-1)^2(2\tau-1)(\tau^3-(\tau-1)^3)\times$$
$$(-\tau^4(\tau+1)^2((\tau+1)^3+2(2\tau+1)^3)-$$
$$(-2\tau^5(\tau+1)^2(2\tau+1)((\tau+1)^3-\tau^3))\times$$
$$(\tau^4(\tau-1)^2((\tau-1)^3+2(2\tau-1)^3)=$$
$$4\tau^9(\tau^2-1)^2(102\tau^6-68\tau^4+21\tau^2-3).$$
\begin{equation}\label{eq:7w1132} 
w^{(1)}_{1,3,2}(\nu)=
\end{equation}
$$-\det\left(\matrix
a_{1,0,1,1}^{\ast(1)\ast}(1;-\nu-2)&a_{1,0,1,2}^{\ast(1)\ast}(1;-\nu-2)\\
a_{1,0,1,1}^{\ast1\ast}(1;\nu)&a_{1,0,1,2}^{\ast1\ast}(1;\nu)
\endmatrix\right)=$$ 
$$a_{1,0,1,1}^{\ast\ast}(1;\nu)a_{1,0,1,2}^{\ast\ast}(1;-\nu-2)-
a_{1,0,1,2}^{\ast\ast}(1;\nu)a_{1,0,1,1}^{\ast\ast}(1;-\nu-2)=$$
$$
-2\tau^4(17\tau^3-27\tau^2+15\tau-3)8\tau^3(6\tau^3+9\tau^2+5\tau+1)-
$$
$$
2\tau^4(17\tau^3+27\tau^2+15\tau+3)8\tau^3(6\tau^3-9\tau^3-5\tau^2+1)=
$$
$$-32\tau^7(102\tau^6-68\tau^4+21\tau^2-3),$$
\begin{equation}\label{eq:7w1232} 
w^(1)_{2,3,2}(\nu)=
\end{equation}
$$-\det\left(\matrix
a_{1,0,2,1}^{\ast1\ast}(1;-\nu-2)&a_{1,0,2,2}^{\ast1\ast}(1;-\nu-2)\\
a_{1,0,2,1}^{\ast1\ast}(1;\nu)&a_{1,0,2,2}^{\ast1\ast}(1;\nu)
\endmatrix\right)=$$
$$\det\left(\matrix
a_{1,0,2,2}^{\ast1\ast}(1;-\nu-2)&a_{1,0,2,1}^{\ast1\ast}(1;-\nu-2)\\
a_{1,0,2,2}^{\ast1\ast}(1;\nu)&a_{1,0,2,1}^{\ast1\ast}(1;\nu)
\endmatrix\right)=$$
$$\det\left(\matrix
a_{1,0,2,2}^{\ast1\ast}(1;-\nu-2)&-a_{1,0,2,3}^{\ast\ast}(1;-\nu-2)\\
a_{1,0,2,2}^{\ast\ast}(1;\nu)&-a_{1,0,2,3}^{\ast\ast}(1;\nu)
\endmatrix\right)=-2\tau^{-2}w^{(0)}_{1,3,2}(\nu).$$
In view of (\ref{eq:7.10}), (\ref{eq:7.14}), (\ref{eq:7w0131}),
 (\ref{eq:7w0132}), (\ref{eq:7w0133}),
\begin{equation}\label{eq:7.15}
(\tau+1)^4(2\tau-1)^3\delta f_{1,0,k}(1,\nu+1)-
\end{equation}
$$4\tau(68\tau^6-45\tau^4+12\tau^2-1)\delta f_{1,0,k}(1,\nu)-$$
$$(\tau-1)^4(2\tau+1)^3\delta f_{\alpha,0,k}(1,\nu-1)=0,$$
where $k=1,\,3$ and
$\nu\in M_1^{\ast\ast\ast\ast}=((-\infty,-3]\cup[1,+\infty))\cap{\mathbb Z}.$
According to the equalities (\ref{eq:7.10}), (\ref{eq:7.14}),
 (\ref{eq:7w0231}),
 (\ref{eq:7w0232}), (\ref{eq:7w0233}),
\begin{equation}\label{eq:7.16}
(\tau+1)^3(2\tau-1)(\tau^3-(\tau-1)^3)\delta^2f_{1,0,k}(1,\nu+1)-
\end{equation}
$$2(102\tau^6-68\tau^4+21\tau^2-3)\delta^2f_{1,0,k}(1,\nu)+$$
$$
(\tau-1)^3(2\tau+1)((\tau+1)^3-\tau^3)\delta^2f_{\alpha,0,k}(1,\nu-1)=0.
$$
In view of (\ref{eq:7.10}), (\ref{eq:7.14}), (\ref{eq:7w1131}),
 (\ref{eq:7w1132}), (\ref{eq:7w1133}),
$$
8\tau^3(2\tau-1)(\tau^3-(\tau-1)^3)(2\tau^2(\tau+1)^5) f_{1,0,k}(1,\nu+1)-
$$
$$32\tau^7(102\tau^6-68\tau^4+21\tau^2-3)f_{1,0,k}(1,\nu)+$$
$$
(-8\tau^3(2\tau+1)((\tau+1)^3-\tau^3))(-2\tau^2(\tau-1)^5f_{1,0,k}(1,\nu-1)=0,
$$
\begin{equation}\label{eq:7.15z}
(\tau+1)^5(2\tau-1)(\tau^3-(\tau-1)^3)f_{1,0,k}(1,\nu+1)-
\end{equation}
$$2\tau^2(102\tau^6-68\tau^4+21\tau^2-3)f_{1,0,k}(1,\nu)+$$
$$(\tau-1)^5(2\tau+1)((\tau+1)^3-\tau^3))f_{1,0,k}(1,\nu-1)=0,$$
where $k=1,\,3$ and
$\nu\in M_1^{\ast\ast\ast\ast}=((-\infty,-3]\cup[1,+\infty))\cap{\mathbb Z}.$

Since
$\delta^2f_{1,0,k}(1,\nu+\varepsilon)=
\frac12(\tau+\varepsilon)^2f_{1,0,k}(1,\nu+\varepsilon))$
if $k=1,\,3,\,\varepsilon^3=\varepsilon$ we can replace 
$\delta^2f_{1,0,k}(1,\nu+\varepsilon)$ by 
$(1/2)(\tau+\varepsilon)^2f_{1,0,k}(1,\nu+\varepsilon))$
in (\ref{eq:7.15}). Then we come to the equality
$$
\frac12(\tau+1)^5(2\tau-1)(\tau^3-(\tau-1)^3)f_{1,0,k}(1,\nu+1)-
$$
$$\tau^2(102\tau^6-68\tau^4+21\tau^2-3)\delta^2f_{1,0,k}(1,\nu)+$$
$$
\frac12(\tau-1)^5(2\tau+1)((\tau+1)^3-\tau^3)\delta^2f_{\alpha,0,k}(1,\nu-1)=0,
$$
this equality is equivalent to the equality (\ref{eq:7.15z}).
The equalities (\ref{eq:7w1231}), (\ref{eq:7w1232}) and (\ref{eq:7w1233})
show that we get in case $\kappa=1$ the same result for the relation
between $\delta^2f_{1,0,k}^\ast(1,\nu+1),\,\delta^2f_{1,0,k}^\ast(1,\nu)$
and $\delta^2f_{1,0,k}(1,\nu-1)$ with $k=1,\,3,$ (relation (\ref{eq:7.16})),
as it have been obtained in the case $\kappa=0.$

So, in the case $\kappa=0$ we obtain the same results as in the
case $\kappa=1.$
Since $$f_{1,0,k}(1,\nu)=\frac1{(\nu+1)^2}f_{1,0,k}^\ast(1,\nu),$$
it follows from (\ref{eq:7.15}) -- ((\ref{eq:7.16}) that
\begin{equation}\label{eq:7.15a}
(\tau+1)^2\tau(2\tau-1)^3\delta f_{1,0,k}^\ast(1,\nu+1)-
\end{equation}
$$4(68\tau^6-45\tau^4+12\tau^2-1)\delta f_{1,0,k}^\ast(1,\nu)+$$
$$(\tau-1)^2\tau(2\tau+1)^3\delta f_{\alpha,0,k}^\ast(1,\nu-1)=0,$$ 
\begin{equation}\label{eq:7.16a}
(\tau+1)\tau^2(2\tau-1)(\tau^3-(\tau-1)^3)\delta^2f_{1,0,k}^\ast(1,\nu+1)-
\end{equation}
$$2(102\tau^6-68\tau^4+21\tau^2-3)\delta^2f_{1,0,k}^\ast(1,\nu)+$$
$$
(\tau-1)\tau^2(2\tau+1)
((\tau+1)^3-\tau^3)\delta^2f_{\alpha,0,k}^\ast(1,\nu-1)=0,
$$
where $k=1,3$
Before to complete the proof of Theorem A, we want to check equalities
(\ref{eq:7.15a}) and (\ref{eq:7.16a}) for $\nu=1,\,k=3.$
If $\nu=1,$ then the left part of the equality (\ref{eq:7.15a})
take the form 
$$
486\delta f_{1,0,k}^\ast(1,2)-14716\delta f_{1,0,k}^\ast(1,1)+
250\delta f_{1,0,k}^\ast(1,0)),
$$
where $k=1,\,3;$
in view of (\ref{eq:a1r1i2n0}) -- (\ref{eq:a1r2i4n2}),
$$\beta^{\ast(1)}_{1,0,2}(1;2)=1029,\,\beta^{\ast(1)}_{1,0,2}(1;1)=34,\,
\beta^{\ast(1)}_{1,0,2}(1;0)=1,$$
$$\beta^{\ast(1)}_{1,0,4}(1;2)=\frac{14843}6,\,
\beta^{\ast(1)}_{1,0,4}(1;1)=\frac{327}4,\,
\beta^{\ast(1)}_{1,0,4}(1;0)=3,$$
$$
486\beta^{\ast(1)}_{1,0,2}(1;2)-14716\beta^{\ast(1)}_{1,0,2}(1;1)+
250\beta^{\ast(1)}_{1,0,2}(1;0)=
$$
$$2(243\times1029-7358\times34+125)=0,$$
$$
486\beta^{\ast(1)}_{1,0,4}(1;2)-14716\beta^{\ast(1)}_{1,0,4}(1;1)+
250\beta^{\ast(1)}_{1,0,4}(1;0)=
$$
$$81\times14843-3679\times327+3\times250=0;$$
therefore, in view of (\ref{eq:15.4.hr1f}), the equality (\ref{eq:7.15a})
holds for $\nu=1.$ 
If $\nu=1,$ then the left part of the equality (\ref{eq:7.16a})
take the form 
$$252\delta^2f_{1,0,k}^\ast(1,2)-
11042\delta^2f_{1,0,k}^\ast(1,\nu)+380\delta^2f_{\alpha,0,k}^\ast(1,0),$$
in view of (\ref{eq:a1r2i2n0}) -- (\ref{eq:a1r2i4n2}),
$$\beta^{\ast(2)}_{1,0,2}(1;2)=2277,\,\beta^{\ast(2)}_{1,0,2}(1;1)=52,\,
\beta^{\ast(2)}_{1,0,2}(1;0)=1,$$
$$\beta^{\ast(2)}_{1,0,4}(1;2)=\frac{32845}6,\,
\beta^{\ast(2)}_{1,0,4}(1;1)=125,\,\beta^{\ast(2)}_{1,0,4}(1;0)=2,$$
$$252\beta^{\ast(2)}_{1,0,2}(1;2)-11042\beta^{\ast(2)}_{1,0,2}(1;1)+
380\beta^{\ast(2)}_{1,0,2}(1;0)=$$
$$252\times2277-11042\times52+380\times1=$$
$$4(63\times2277-11042\times13+95)=0,$$
$$252\beta^{\ast(2)}_{1,0,4}(1;2)-
11042\beta^{\ast(2)}_{1,0,4}(1;1)+
380\beta^{\ast(2)}_{1,0,4}(1;0)=
$$
$$42\times32845-11042\times125+380\times2=$$
$$10(21\times6569-5521\times25+76)=0;$$
therefore, in view of (\ref{eq:15.4.hr1f}), the equality (\ref{eq:7.16a})
holds for $\nu=1.$
Let us consider the equations 
\begin{equation}\label{eq:7.17}
(\tau+1)^2\tau(2\tau-1)^3x_{nu+1}-
\end{equation}
$$4(68\tau^6-45\tau^4+12\tau^2-1)x_\nu+$$
$$(\tau-1)^2\tau(2\tau+1)^3x_{\nu-1}=0,$$
\begin{equation}\label{eq:7.18}
(\tau+1)(2\tau-1)\tau^2(\tau^3-(\tau-1)^3)x_{nu+1}-
\end{equation}
$$2(102\tau^6-68\tau^4+21\tau^2-3)x_nu+$$
$$
(\tau-1)(2\tau+1)\tau^2((\tau+1)^3-\tau^3)x_{nu-1}=0,
$$
where $\nu\in{\mathbb N}_0,\,\tau=\nu+1.$
 It follows from (\ref{eq:7.15a}) that $x_\nu=\delta f_{1,0,k}^\ast(1,\nu)$
satisfies to the equation (\ref{eq:7.17}) for $\nu\in{\mathbb N}$ and
 fixed $k\in\{1,\,3\}.$
 It follows from (\ref{eq:7.16a}) that $x_\nu=\delta^2 f_{1,0,k}^\ast(1,\nu)$
satisfies (\ref{eq:7.17}) for $\nu\in{\mathbb N}$ and
 fixed $k\in\{1,\,3\}.$
 Both equations (\ref{eq:7.17}) and (\ref{eq:7.18}) are difference
 equations of Poincar\'e type with characteristic polynomial
 $\lambda^2-34\lambda+1.$
 Hence, if $\{x_\nu\}_{\nu=1}^{+\infty}$ is a non-zero solution 
  of some of these equations, $\varepsilon\in(0,1),$
   then there are $C_1(\varepsilon)>0$ and $C_2(\varepsilon)>0$ such that
 only two possibilities exist:
 \begin{equation}\label{eq:7.bf1}
\frac{C_1(\varepsilon)}{\left(1+\sqrt{2}\right)^{4\nu(1+\varepsilon)}}\le
\vert x_\nu\vert\le
\frac{C_2(\varepsilon)}{\left(1+\sqrt{2}\right)^{4\nu(1-\varepsilon)}}
\end{equation}
for all $\nu\in{\mathbb N}$ or
 \begin{equation}\label{eq:7.bf2}
C_1(\varepsilon)\left(1+\sqrt{2}\right)^{4\nu(1-\varepsilon)}\le
\vert x_\nu\vert\le
C_2(\varepsilon)\left(1+\sqrt{2}\right)^{4\nu(1+\varepsilon)}.
\end{equation}
for all $\nu\in{\mathbb N}.$
In view of (\ref{eq:15.4.hr8}), if $x_\nu=\delta^r f_{1,0,1}^\ast(1,\nu)$
with $r=1,2,$ then (\ref{eq:7.bf1}) is impossible.
In view of (\ref{eq:d}) with $\alpha\in{\mathbb N},$
(\ref{eq:15.4.hr}) and (\ref{eq:15.4.hr1f}) with $j=1,$
$$\delta^r f_{1,0,3}^\ast(1,\nu)=(\nu+1)^2O(1);$$ hence,
 if $x_\nu=\delta^r f_{1,0,3}^\ast(1,\nu)$
with $r=1,2,$ then (\ref{eq:7.bf2}) is impossible.
In view of (\ref{eq:15.4.hr6}) with $j=1$ and (\ref{eq:15.4.hr8}),
$x_\nu=\beta^{\ast(1)}_{1,0,4}(1,\nu)$ and
$x_\nu=\beta^{\ast(2)}_{1,0,4}(1,\nu)$ are
solutions for $\nu\in{\mathbb N}$
of equations respectively (\ref{eq:7.17}) and (\ref{eq:7.18}), moreover
(\ref{eq:7.bf2}) take place for these solutions.
Hence,
 \begin{equation}\label{eq:7.bf3}
\frac{C_1(\varepsilon)/C_2(\varepsilon)}
{\left(1+\sqrt{2}\right)^{8\nu(1+\varepsilon)}}\le 
\bigg\vert2\zeta(3)-\frac{\beta^{\ast(r)}_{1,0,4}(1,\nu)}
{\beta^{\ast(r)}_{1,0,2}(1,\nu)}\bigg\vert\le
\end{equation}
$$\frac{C_2(\varepsilon)/C_1(\varepsilon)}
{\left(1+\sqrt{2}\right)^{8\nu(1-\varepsilon)}}$$
for $r=1,\,2.$
The equations (\ref{eq:7.17}) and (\ref{eq:7.18}) are equivalent to
the equations respectively
\begin{equation}\label{eq:0c001}
x_{\nu+1}=b^{(1)}_{\nu+1}x_\nu+a^{(1)}_{\nu+1}x_{\nu-1},
\end{equation}
and
\begin{equation}\label{eq:0c002}
x_{\nu+1}=b^{(2)}_{\nu+1}x_\nu+a^{(2)}_{\nu+1}x_{\nu-1},
\end{equation}
with
 \begin{equation}\label{eq:7.bf4}
b^{(1)}_{\nu+1}=b^{(1)}_\tau=\frac{4(68\tau^6-45\tau^4+12\tau^2-1)}
{(\tau+1)^2\tau(2\tau-1)^3},
\end{equation}
 \begin{equation}\label{eq:7.bf5}
a^{(1)}_{\nu+1}=a^{(1)}_\tau=-\frac{(\tau-1)^2\tau(2\tau+1)^3}
{(\tau+1)^2\tau(2\tau-1)^3},
\end{equation}
 \begin{equation}\label{eq:7.bg4}
b^{(2)}_{\nu+1}=b^{(2)}_\tau=
\end{equation}
$$\frac{2(102\tau^6-68\tau^4+21\tau^2-3)}
{(\tau+1)(2\tau-1)\tau^2(\tau^3-(\tau-1)^3)},$$
 \begin{equation}\label{eq:7.bg5}
a^{(2)}_{\nu+1}=a^{(2)}_\tau=
\end{equation}
$$-\frac{(\tau-1)(2\tau+1)\tau^2((\tau+1)^3-\tau^3)}
{(\tau+1)(2\tau-1)\tau^2(\tau^3-(\tau-1)^3)}$$
and $\nu\in{\mathbb N}.$ To consider the case $\nu=0$ also, we let
 \begin{equation}\label{eq:7.bf45n0}
a^{(1)}_1=\frac{-81}4,\,b^{(1)}_0=\beta^{\ast(1)}_{1,0,4}(1,1)=3,\,b^{(1)}_1=
34=
\end{equation}
$$\frac{4(68\tau^6-45\tau^4+12\tau^2-1)}
{(\tau+1)^2\tau(2\tau-1)^3}\bigg\vert_{\tau=1},$$
 \begin{equation}\label{eq:7.bg45n0}
a^{(2)}_1=21,\,b^{(2)}_0=\beta^{\ast(2)}_{1,0,4}(1,0)=2,\,b^{(2)}_1=52=
\end{equation}
$$\frac{2(102\tau^6-68\tau^4+21\tau^2-3)}
{(\tau+1)(2\tau-1)\tau^2(\tau^3-(\tau-1)^3)}\bigg\vert_{\tau=1}.$$
The equation (\ref{eq:0c001}) is the eqution (\ref{eq:0c00})
with $a^{(1)}_\nu$ and $b^{(1)}_\nu$ in the role of respectively
$a_\nu$ and $b_\nu.$
The equation (\ref{eq:0c002}) is the eqution (\ref{eq:0c00})
with $a^{(2)}_\nu$ and $b^{(2)}_\nu$ in the role of respectively
$a_\nu$ and $b_\nu.$
Let $x_\nu=Q^{(1)}_{1,\nu}$ and $x_\nu=Q^{(2)}_\nu,$
 where $\nu\in{\mathbb N}_0,$ are solutions 
of the equations respectively (\ref{eq:0c001}) and (\ref{eq:0c002})
  with initial values respectively
   $x_0=1,\,x_1=b^{(1)}_1=34$ and $x_0=1,\,x_1=b^{(2)}_1=52.$
 Then we see that $x_\nu=\beta^{\ast(r)}_{1,0,2}(1,\nu)$ 
and $x_\nu=Q^{(r)}_\nu$ with fixed $r\in\{1,\,2\}$ are solutions
of the same equation with the same initial values. Therefore
\begin{equation}\label{eq:7.bfg6}
Q^{(r)}_\nu=\beta^{\ast(r)}_{1,0,2}(1,\nu)
\end{equation}
for $r=1,\,2$ and $\nu\in{\mathbb N}_0.$
Let $x_\nu=P^{(1)}_\nu$ and $x_\nu=P^{(2)}_\nu,$
 where $\nu\in{\mathbb N}_0,$ are solutions
of the equations respectively (\ref{eq:0c001}) and (\ref{eq:0c002})
  with initial values respectively
$$x_0=b^{(1)}_0,\,x_1=b^{(1)}_1b^{(1)}_0+a^{(1)}_1=\frac{327}4=
\beta^{\ast(1)}_{1,0,4}(1,1)$$
and
$$x_0=b^{(2)}_0,\,x_1=b^{(2)}_1b^{(2)}_0+a^{(2)}_1=125=
\beta^{\ast(2)}_{1,0,4}(1,1).$$
Since $x_\nu=\beta^{\ast(r)}_{1,0,4}(1,\nu)$ 
and $x_\nu=P^{(r)}_\nu$ with fixed $r\in\{1,\,2\}$ are solutions
of the same equation with the same initial values, it follows that
\begin{equation}\label{eq:7.bfg7}
P^{(r)}_\nu=\beta^{\ast(r)}_{1,0,4}(1,\nu)
\end{equation}
for $r=1,\,2$ and $\nu\in{\mathbb N}_0.$
 In view of (\ref{eq:7.bf3}) with $r=1,$
 \begin{equation}\label{eq:7.bf3a}
\frac{C_1(\varepsilon)/C_2(\varepsilon)}
{\left(1+\sqrt{2}\right)^{8\nu(1+\varepsilon)}} 
\le\bigg\vert2\zeta(3)-\frac{P^{(r)}_\nu}
{Q^{(r)}_\nu}\bigg\vert\le
\frac{C_2(\varepsilon)/C_1(\varepsilon)}
{\left(1+\sqrt{2}\right)^{8\nu(1-\varepsilon)}}.
\end{equation}
We put in the equation (\ref{eq:0c00}) $\xi_\nu=x_\nu/d_\nu$
with some $d_\nu\ne0$ for $\nu\in{\mathbb N}_0,$ and
we put $d_0=d_{-1}=1.$
Then we obtain equation
\begin{equation}\label{eq:0c003}
\xi_{\nu+1}=\beta^\vee_{\nu+1}\xi_\nu+
\alpha_{\nu+1}\xi_{\nu-1}
\end{equation}
with $\xi_0=x_0=b_0,$
\begin{equation}\label{eq:7.bh7}
\alpha_{\nu+1}=a_{\nu+1}\frac{d_{\nu+1}}{d_{\nu-1}},
\beta_{\nu+1}=b_{\nu+1}\frac{d_{\nu+1}}{d_\nu}
\end{equation}
where $\nu\in{\mathbb N}_0.$
Clearly, continuous fraction connected with the equation (\ref{eq:0c003})
has the same convergents as
continuous fraction corresponding to the equation (\ref{eq:0c00}).
We apply this transformatin for the equations (\ref{eq:0c001})
and (\ref{eq:0c002}) 
instead of equation (\ref{eq:0c00}) with respectively
$$d_\nu=d^{(1)}_\nu=\prod\limits_{k=1}^\nu k(k+1)^2(2k-1)^3,$$
$$
d_\nu=d^{(2)}_\nu=\prod\limits_{k=1}^\nu (k+1)(2k-1)k^2(k^3-(k-1)^3),
$$
Then we have
$$\frac{d^{(1)}_{\nu+1}}{d^{(1)}_\nu}=(\nu+1)(\nu+2)^2(2\nu+1)^3=
\tau(\tau+1)^2(2\tau-1)^3,$$
$$\frac{d^{(2)}_{\nu+1}}{d^{(2)}_\nu}=
(\nu+2)(2\nu+1)(\nu+1)^2((\nu+1)^3-\nu^3)=$$
$$(\tau+1)(2\tau-1)\tau^2(\tau^3-(\tau-1)^3),$$
if $\nu\in{\mathbb N}_0,$
$$\frac{d^{(1)}_{\nu+1}}{d^{(1)}_{\nu-1}}=
\tau(\tau+1)^2(2\tau-1)^3\times
(\tau-1)\tau^2(2\tau-3)^3,$$
$$\frac{d^{(2)}_{\nu+1}}{d^{(2)}_{\nu-1}}=
(\tau+1)(2\tau-1)\tau^2(\tau^3-(\tau-1)^3)\times$$
$$\tau(2\tau-3)(\tau-1)^2((\tau-1)^3-(\tau-2)^3),$$
if $\nu\in{\mathbb N},$
$$\frac{d^{(1)}_1}{d^{(1)}_{-1}}=4,
\frac{d^{(2)}_1}{d^{(2)}_{-1}}=2,$$
and we obtain two expansions of the number $2\zeta(3)$
in continuous fraction:
one expansion with $a_1=-81,\, b_0=3,$ 
$$a_{\nu+1}=a_\tau=
-\frac{(\tau-1)^2\tau(2\tau+1)^3}
{(\tau+1)^2\tau(2\tau-1)^3}\times$$
$$\tau(\tau+1)^2(2\tau-1)^3\times
(\tau-1)\tau^2(2\tau-3)^3=$$
$$-(\tau(\tau+1)(4\tau^2-4\tau-3))^3$$
for $\nu\in{\mathbb N},$
$$ b_{\nu+1}=b_\tau=
4(68\tau^6-45\tau^4+12\tau^2-1)$$
for $\nu\in{\mathbb N}_0,$
and another expansion with 
$$a_1=42,\, b_0=2,\,a_{\nu+1}=a_\tau=$$
$$-\frac{(\tau-1)(2\tau+1)\tau^2((\tau+1)^3-\tau^3)}
{(\tau+1)(2\tau-1)\tau^2(\tau^3-(\tau-1)^3)}
\times$$
$$(\tau+1)(2\tau-1)\tau^2(\tau^3-(\tau-1)^3)\times$$
$$\tau(2\tau-3)(\tau-1)^2((\tau-1)^3-(\tau-2)^3)=$$
$$
-(\tau-1)^3\tau^3(4\tau^2-4\tau-3)\times$$
$$((\tau+1)^3-\tau^3)((\tau-1)^3-(\tau-2)^3)$$
for $\nu\in{\mathbb N},$
$$ b_{\nu+1}=b_\tau=
2(102\tau^6-68\tau^4+21\tau^2-3)$$
for $\nu\in{\mathbb N}_0,$
$\blacksquare$
{\begin{center}\large\bf References.\end{center}}
\footnotesize
\vskip4pt
\refstepcounter{r}\noindent[\ther]
R.Ap\'ery, Interpolation des fractions continues\\
\hspace*{3cm}  et irrationalite de certaines constantes,\\
\hspace*{3cm} Bulletin de la section des sciences du C.T.H., 1981, No 3,
 37 -- 53;
\label{r:cd}\\
\refstepcounter{r}
Oskar Perron. Die Lehre von den Kettenbr\"uche.
Dritte, verbesserte und erweiterte Auflage. 1954
B.G.Teubner Verlaggesellshaft. Stuttgart. 
\label{r:ce}\\
\refstepcounter{r}
Yu.V. Nesterenko. A Few Remarks on $\zeta(3),$\\
\hspace*{3cm}Mathematical Notes, Vol 59, No 6, 1996,\\
\hspace*{3cm} Matematicheskie Zametki, Vol 59, No 6, pp. 865 -- 880 1996\\
\hspace*{3cm} (in Russian).
\label{r:cf}\\
\refstepcounter{r}
\noindent [\ther] L.A.Gutnik, On linear forms with coefficients
in ${\mathbb N}\zeta(1+{\mathbb N})$\\
\hspace*{3cm} (the detailed version,part 3),
 Max-Plank-Institut f\"ur Mathematik,\\
\hspace*{3cm} Bonn, Preprint Series, 2002, 57, 1 -- 33;
\label{r:bh}\\
\refstepcounter{r}
\noindent [\ther] \rule{2cm}{.3pt}, On the measure of nondiscreteness
of some modules,\\
\hspace*{3cm} Max-Plank-Institut f\"ur Mathematik, Bonn,\\
\hspace*{3cm} Preprint Series, 2005, 32, 1 -- 51.
*\label{r:caj2}\\
\refstepcounter{r}
\noindent [\ther] \rule{2cm}{.3pt}, On the Diophantine approximations\\
\hspace*{3cm} of logarithms in cylotomic fields.\\
\hspace*{3cm} Max-Plank-Institut f\"ur Mathematik, Bonn,\\
\hspace*{3cm} Preprint Series, 2006, 147, 1 -- 36.
\label{r:daj}\\
\refstepcounter{r}
\noindent [\ther] \rule{2cm}{.3pt}, On some systems of difference
equations. Part 1.\\
\hspace*{3cm} Max-Plank-Institut f\"ur Mathematik, Bonn,\\
\hspace*{3cm} Preprint Series, 2006, 23, 1 -- 37.
\label{r:caj3}\\
\refstepcounter{r}
\noindent [\ther] \rule{2cm}{.3pt}, On some systems of difference
equations. Part 2.\\
\hspace*{3cm} Max-Plank-Institut f\"ur Mathematik, Bonn,\\
\hspace*{3cm} Preprint Series, 2006, 49, 1 -- 31.
\label{r:caj4}\\
\refstepcounter{r}
\noindent [\ther] \rule{2cm}{.3pt}, On some systems of difference
equations. Part 3.\\
\hspace*{3cm} Max-Plank-Institut f\"ur Mathematik, Bonn,\\
\hspace*{3cm} Preprint Series, 2006, 91, 1 -- 52.
\label{r:caj5}\\
\refstepcounter{r}
\noindent [\ther] \rule{2cm}{.3pt}, On some systems of difference
equations. Part 4.\\
\hspace*{3cm} Max-Plank-Institut f\"ur Mathematik, Bonn,\\
\hspace*{3cm} Preprint Series, 2006, 101, 1 -- 49.
\label{r:caj6}\\
\refstepcounter{r}
\noindent [\ther] \rule{2cm}{.3pt}, On some systems of difference
equations. Part 5.\\
\hspace*{3cm} Max-Plank-Institut f\"ur Mathematik, Bonn,\\
\hspace*{3cm} Preprint Series, 2006, 115, 1 -- 9.
\label{r:caj7}\\
\refstepcounter{r}
\noindent [\ther] \rule{2cm}{.3pt}, On some systems of difference
equations. Part 6.\\
\hspace*{3cm} Max-Plank-Institut f\"ur Mathematik, Bonn,\\
\hspace*{3cm} Preprint Series, 2007, 16, 1 -- 30.
\label{r:caj8}\\
\refstepcounter{r}
\noindent [\ther] \rule{2cm}{.3pt}, On some systems of difference
equations. Part 7.\\
\hspace*{3cm} Max-Plank-Institut f\"ur Mathematik, Bonn,\\
\hspace*{3cm} Preprint Series, 2007, 53, 1 -- 40.
\label{r:caj9}\\
\refstepcounter{r}
\noindent [\ther] \rule{2cm}{.3pt}, On some systems of difference
equations. Part 8.\\
\hspace*{3cm} Max-Plank-Institut f\"ur Mathematik, Bonn,\\
\hspace*{3cm} Preprint Series, 2007, 64, 1 -- 44.
\label{r:caj10}\\
\refstepcounter{r}
\noindent [\ther] \rule{2cm}{.3pt}, On some systems of difference
equations. Part 9.\\
\hspace*{3cm} Max-Plank-Institut f\"ur Mathematik, Bonn,\\
\hspace*{3cm} Preprint Series, 2007, 129, 1 -- 36.
\label{r:caj11}\\
\refstepcounter{r}
\noindent [\ther] \rule{2cm}{.3pt}, On some systems of difference
equations. Part 10.\\
\hspace*{3cm} Max-Plank-Institut f\"ur Mathematik, Bonn,\\
\hspace*{3cm} Preprint Series, 2007, 131 , 1 -- 33. 
\label{r:caj12}\\
\refstepcounter{r}
\noindent [\ther] \rule{2cm}{.3pt}, On some systems of difference
equations. Part 11.\\
\hspace*{3cm} Max-Plank-Institut f\"ur Mathematik, Bonn,\\
\hspace*{3cm} Preprint Series, 2008, 38 , 1 -- 45. 
\label{r:caj13}
\vskip 10pt
 {\it E-mail:}{\sl\ gutnik$@@$gutnik.mccme.ru}
\end{document}